\documentclass[11pt]{article}

\usepackage{amscd,amsmath, amssymb, fancyhdr, epsfig, color, tocloft, mathtools}
\usepackage{dutchcal}

\usepackage[backref=page]{hyperref}
\renewcommand*{\backrefalt}[4]{%
    \ifcase #1 (Not cited.)%
    \or        (Cited on page~#2.)%
    \else      (Cited on pages~#2.)%
    \fi}

\hypersetup{
     colorlinks   = true,
     citecolor    = magenta,
     linkcolor    = blue,
     urlcolor     = magenta	
}

\newcommand{\version}{version 2.2.1,\ \ September 27, 2021}
\makeatletter
\def\x@arrow{\DOTSB\Relbar}
\def\xlongequalsignfill@{\arrowfill@\x@arrow\Relbar\x@arrow}
\providecommand{\xlongequal}[2][]{%
	\ext@arrow 0099\xlongequalsignfill@{#1}{#2}}
\def\xlongrightarrowfill@{\arrowfill@\relbar\relbar\longrightarrow}
\newcommand{\xlongrightarrow}[2][]{%
        \ext@arrow 0099\xlongrightarrowfill@{#1}{#2}}
\makeatother

\numberwithin{equation}{section}

\def\eqref#1{(\ref{#1})}
\newcommand{\goth}{\mathfrak}
\newcommand{\g}{{\mathfrak g}}
\renewcommand{\a}{{\goth a}}

\newcommand{\Z}{{\mathbb Z}}
\newcommand{\C}{{\mathbb C}}
\newcommand{\R}{{\mathbb R}}

\def\1{\sqrt{-1}\:}
\newcommand{\restrict}[1]{{\left|_{{\phantom{|}\!\!}_{#1}}\right.}}
\newcommand{\cntrct}                
{\hspace{2pt}\raisebox{1pt}{\text{$\lrcorner$}}\hspace{2pt}}
\newcommand{\arrow}{{\:\longrightarrow\:}}

\renewcommand{\phi}{\varphi}
\renewcommand{\epsilon}{\varepsilon}
\renewcommand{\geq}{\geqslant}
\renewcommand{\leq}{\leqslant}

\newcommand{\ev}{{\rm even}}
\newcommand{\odd}{{\rm odd}}

\newcommand{\im}{\operatorname{im}}
\newcommand{\End}{\operatorname{End}}

\newcommand{\Id}{\operatorname{Id}}

\newcommand{\kah}{{\operatorname{\text{\sf kah}}}}
\newcommand{\sas}{{\operatorname{\text{\sf sas}}}}
\newcommand{\bas}{{\operatorname{\text{\sf bas}}}}
\newcommand{\hor}{{\operatorname{\text{\sf hor}}}}
\renewcommand{\vert}{{\operatorname{\text{\sf vert}}}}

\newcommand{\Iso}{\operatorname{Iso}}

\newcommand{\Spin}{\operatorname{Spin}}
\newcommand{\Av}{\operatorname{Av}}

\newcommand{\Diff}{\operatorname{Diff}}

\newcommand{\rk}{\operatorname{rk}}
\newcommand{\Lie}{\operatorname{Lie}}
\newcommand{\Sym}{\operatorname{Sym}}

\newcommand{\GL}{\operatorname{GL}}

\newcounter{Mycounter}[section]
\newcounter{lemma}[section]
\renewcommand{\thelemma}{{Lemma \thesection.\arabic{lemma}}}
\newcommand{\lemma}{%
     \setcounter{lemma}{\value{Mycounter}}
     \refstepcounter{lemma}
     \stepcounter{Mycounter}
     {\noindent \bf \thelemma:\ }}

\newcounter{claim}[section]
\renewcommand{\theclaim}{{Claim \thesection.\arabic{claim}}}
\newcommand{\claim}{%
     \setcounter{claim}{\value{Mycounter}}
     \refstepcounter{claim}
     \stepcounter{Mycounter}
     {\noindent \bf \theclaim:\ }}

\newcounter{sublemma}[section]
\newcounter{corollary}[section]
\renewcommand{\thecorollary}{{Corollary \thesection.\arabic{corollary}}}
\newcommand{\corollary}{%
     \setcounter{corollary}{\value{Mycounter}}
     \refstepcounter{corollary}
     \stepcounter{Mycounter}
     {\noindent \bf \thecorollary:\ }}

\newcounter{theorem}[section]
\renewcommand{\thetheorem}{{Theorem \thesection.\arabic{theorem}}}
\newcommand{\theorem}{%
     \setcounter{theorem}{\value{Mycounter}}
     \refstepcounter{theorem}
     \stepcounter{Mycounter}
     {\noindent \bf \thetheorem:\ }}

\newcounter{conjecture}[section]
\newcounter{proposition}[section]
\renewcommand{\theproposition} {{Proposition \thesection.\arabic{proposition}}}
\newcommand{\proposition}{%
     \setcounter{proposition}{\value{Mycounter}}
     \refstepcounter{proposition}
     \stepcounter{Mycounter}
     {\noindent \bf \theproposition:\ }}

\newcounter{definition}[section]
\renewcommand{\thedefinition} {{Definition~\thesection.\arabic{definition}}}
\newcommand{\definition}{%
     \setcounter{definition}{\value{Mycounter}}
     \refstepcounter{definition}
     \stepcounter{Mycounter}
     {\noindent \bf \thedefinition:\ }}

\newcounter{example}[section]
\renewcommand{\theexample}{{Example \thesection.\arabic{example}}}
\newcommand{\example}{%
     \setcounter{example}{\value{Mycounter}}
     \refstepcounter{example}
     \stepcounter{Mycounter}
     {\noindent \bf \theexample:\ }}

\newcounter{remark}[section]
\renewcommand{\theremark}{{Remark \thesection.\arabic{remark}}}
\newcommand{\remark}{%
     \setcounter{remark}{\value{Mycounter}}
     \refstepcounter{remark}
     \stepcounter{Mycounter}
     {\noindent \bf \theremark:\ }}

\newcounter{problem}[section]
\newcounter{question}[section]
\makeatletter

\@addtoreset{equation}{section}
\@addtoreset{footnote}{section}
\makeatother

\def\blacksquare{\hbox{\vrule width 5pt height 5pt depth 0pt}}
\def\endproof{\blacksquare}

\newcommand{\proof}{{\bf Proof: \ }}
\newcommand{\pstep}{{\bf Proof. Step 1: \ }}

\begin{document}

\begin{center}
{\Large\bf  Supersymmetry and Hodge theory on\\ Sasakian and Vaisman manifolds}\\[5mm]
{\large
Liviu Ornea\footnote{Liviu Ornea is  partially supported by a grant of Ministry of Research and Innovation, CNCS - UEFISCDI, project number PN-III-P4-ID-PCE-2016-0065, within PNCDI III.},  
Misha
Verbitsky\footnote{Misha Verbitsky is partially supported by
by the HSE University Basic Research Program, FAPERJ E-26/202.912/2018 
and CNPq - Process 313608/2017-2.\\[1mm]

}\\[4mm]

}

\end{center}

{\small
\hspace{0.15\linewidth}
\begin{minipage}[t]{0.7\linewidth}
{\bf Abstract} \\ 
Sasakian manifolds are odd-dimensional counterpart to
K\"ahler manifolds. They can be defined as contact
manifolds equipped with an invariant K\"ahler structure on
their symplectic cone. The quotient of this cone by the
homothety action is a complex manifold called Vaisman.  We
study harmonic forms and Hodge decomposition on  Vaisman
and Sasakian manifolds. We construct a Lie superalgebra
associated to a Sasakian manifold in the same way as the
K\"ahler supersymmetry algebra is associated to a K\"ahler
manifold. We use this construction to produce a
self-contained, coordinate-free proof of the results by
Tachibana, Kashiwada and Sato on the decomposition of
harmonic forms and cohomology of Sasakian and Vaisman
manifolds.  In the last section, we compute the
supersymmetry algebra of Sasakian manifolds explicitly.
\end{minipage}
}

\tableofcontents


\section{Introduction}
\label{_Intro_Section_}


Sasakian manifolds are odd-dimensional
counterparts to  K\"ahler manifolds. They can be
defined as contact manifolds equipped with an invariant
K\"ahler structure on their symplectic cone.
Taking a quotient of this cone by the
homothety action one obtains a complex
manifold called Vaisman. Here we develop a version of
Hodge theory on Sasakian and Vaisman manifolds.

The modern approach to Hodge theory is inspired by the supersymmetry.
Given a manifold with a geometric structure (such as K\"ahler,
hyperk\"ahler,  HKT, $G_2$-manifold and so on),
one takes a bunch of natural operators on the de Rham algebra
(such as the de Rham differential, the Lefschetz ${\goth {sl}}(2)$-triple
etc.) and proves that these operators generate
a finite-dimensional Lie superalgebra. In most cases,
the Laplacian is central in this superalgebra.
This gives an interesting geometric action on the
cohomology.

In this paper we write the natural superalgebra
for a Sasakian manifold. Hodge theory for Sasakian
manifolds is well developed, but most proofs are written
in classical differential geometry style, mixing
coordinate computations  with invariant ones. We wanted to
give a more conceptual proof. 

In the end we arrived at a very simple approach
to the Hodge theory on Sasakian and Vaisman manifolds
(\ref{_harmo_Sasa_decompo_Theorem_}, 
\ref{_Vaisman_harmonic_forms_Theorem_}). 
The Sasakian supersymmetry algebra is
not used in this development, but it is  
interesting in itself, 
and the relations are quite surprising.

When the Sasakian manifold is also Einstein,
the superalgebra structure seems to be simpler.
Using this approach, J. Schmude (\cite{_Schmude_})
obtained a closed formula for the de Rham Laplacian
in terms of the transversal Laplacian operator.

\subsection{Supersymmetry in K\"ahler and non-K\"ahler geometry}

The connection between de Rham calculus on manifolds
with special geometry (such as K\"ahler and hyperk\"ahler)
and their supersymmetry appeared as early as in 1997 
(\cite{_FKS_}). It is well known that extra supersymmetries
of the $\sigma$-model force the target space to acquire
extra geometric structures: the $N=1$ supersymmetry 
implies K\"ahler structure on the target space,
the $N=2$ supersymmetry makes it hyperk\"ahler, and so on.
In \cite{_FKS_}, this supersymmetry was interpreted
in terms of the de Rham calculus on the target space.

In \cite{_Kim_Saberi_}, the connection between 
the supersymmetry and rational homotopy theory is further
expounded, with the constructions of rational homotopy
theory (such as Sullivan's minimal models) interpreted
in terms of quantum mechanics.

For Sasakian manifolds, this approach was pioneered
by \cite{_Tievski_}, who used the transversal K\"ahler relations
to obtain results about rational homotopy of
Sasakian manifolds. This work was
applied to homotopy formality of Sasakian manifolds
(\cite{_BFM_})
and in applications to the geometry of Sasakian
nilmanifolds (\cite{_CDMY_}).

Traditionally, the K\"ahler identities were
obtained using the Levi-Civita connection.
Alternatively, one may show that a K\"ahler
manifold can be approximated up to second order
by a flat space. Using supersymmetry to prove
the K\"ahler identities (Section \ref{_SUSY_Section_}) has many advantages
over either of these approaches. In \cite{_Verbitsky:HKT_}, 
supersymmetry was used to develop the Hodge theory for the HKT 
manifolds (hyperk\"ahler manifolds with torsion),
where the structure tensors are neither preserved by 
the Levi-Civita connection nor admit a second
order approximation. The supersymmetry approach 
was later used in supersymmetric $\sigma$-models
associated to HKT manifolds 
(\cite{_Smilga:supercharges_,_FS:bi_,_FIS:generic_}).

In \cite{_Verbitsky:NK_Hodge_}, the supersymmetry approach was used
to obtain the Hodge decomposition on nearly K\"ahler
manifolds, where the structure tensors are also non-parallel.

In another direction, one may use the supersymmetry
to obtain the superalgebra action on manifolds with
parallel differential forms (\cite{_Verbitsky:G2_forma_}).
The Hodge-theoretic results obtained in this direction
were used in \cite{_Huang:L2_} to study the geometry of complete
$G_2$- and $\Spin(7)$-manifolds with 
$d(\text{linear})$ structure form.

Summing up, finding a superalgebra associated
with a given geometric structure seems to be worthwhile,
for mathematics as well as for mathematical physics.

\subsection{Structure of this paper}

\begin{itemize}
\item In Section \ref{_SUSY_Section_} we relate the Lie superalgebra approach
to de Rham calculus and obtain the K\"ahler identities
which lie at the foundation of  Hodge theory.

\item Section \ref{_Sasak_Section_} introduces the Sasakian manifolds and
explains the basic notions of Sasakian geometry.

\item In Section \ref{_Hattori_Section_} we explain how one employs de Rham calculus
to obtain the Leray-Serre spectral sequence; this approach
was used by A. Hattori in 1960. We need a version of
the Leray-Serre spectral sequence which can be applied
to smooth foliations with fibers which are not necessarily
closed or compact. We compute the differentials
of this spectral sequence for Sasakian
manifolds explicitly, for later use.

\item In Section \ref{_Tra_Kah_Section_} we deal with transversally K\"ahler
structures on smooth foliations. We give simple proofs
of the standard results on transversal (basic) cohomology, 
due to El Kacimi-Alaoui (\cite{_Kacimi_}), much simplified because
we need them only for Vaisman and Sasakian manifolds.

\item In Section \ref{_Hodge_Sasakian_Section_}, we prove the standard results on 
cohomology for Sasakian manifolds, in Section \ref{_Vaisman_Section_}
we introduce the Vaisman manifolds, and in Section \ref{_Hodge_Vaisman_Section_}
we prove the standard results on cohomology for Vaisman manifolds. 
The cohomology calculations for Sasakian and Vaisman manifolds are
very similar and rely on the same homological algebra 
argument. Throughout Section \ref{_Hodge_Vaisman_Section_}, we use the
1-dimensional transversally Sasakian foliation,
generated by the Lee field, and apply the 
superalgebra computations done in Section 
\ref{_Hodge_Sasakian_Section_}
on Sasakian manifolds.
However, the results on Vaisman manifolds
cannot be deduced from the results on Sasakian manifolds,
because the Lee foliation might have non-closed leaves.

\item We finish the paper with an explicit
calculation of the Sasakian superalgebra in Section 
\ref{_SUSY_Sasakian_Section_}. This section depends only
on Sections \ref{_SUSY_Section_}-\ref{_Hattori_Section_}.
\end{itemize}


\section{Lie superalgebras acting on the de Rham algebra}
\label{_SUSY_Section_}


\newcommand{\vecr}{{\vec r}}

\subsection{Lie superalgebras and superderivations}

In the following, all vector spaces and algebras are considered over $\R$. Let $A$ be a $\Z/2\Z$-graded vector space,
\[ A = A^\ev \oplus A^\odd.\] We say that $a\in A$
is {\bf pure} if $a$ belongs to $A^\ev$ or $A^\odd$.
For a pure element $a\in A$, we
write $\tilde a =0$ if $a\in A^\ev$,
and $\tilde a =1$ if $a\in A^\odd$.
Consider a bilinear operator
\[ \{\cdot, \cdot\}:\; A \times A \arrow A,\]
called {\bf supercommutator}. Assume that 
$\{\cdot, \cdot\}$ is {\bf graded an\-ti-\-com\-mu\-ta\-ti\-ve},
that is, satisfies
\[ \{ a, b\} = - (-1)^{\tilde a \tilde b}  \{ b,a \}
\]
for pure $a, b \in A$. Assume, moreover, that
$\{\cdot, \cdot\}$ is compatible with the grading:
the commutator $\{ a, b\}$ is even when both $a$, $b$ are even 
or odd, and odd if one of these elements is odd and another 
is even. We say that $(A, \{\cdot, \cdot\})$ is a {\bf Lie superalgebra} 
if the following identity (called {\bf the graded Jacobi
identity}, or {\bf super Jacobi identity}) 
holds, for all pure elements $a, b, c\in A$: 
\begin{equation}
\{ a, \{ b, c\}\} = \{ \{ a, b\}, c\} + (-1)^{\tilde a \tilde b}\{b, \{ a,c \}\}.
\end{equation}
Up to a sign, this is the usual Jacobi identity.

Every reasonable property of Lie algebras has a 
natural analogue for Lie superalgebras, using the following  rule of thumb: 
every time one would exchange two elements $a$ and $b$,
one adds a multiplier $(-1)^{\tilde a \tilde b}$.

\hfill

\example 
Let $V= V^\ev\oplus V^\odd$ be a $\Z/2\Z$-graded vector space,
and $\End (V)$ its space of endomorphisms, equipped with the induced grading.
We define a supercommutator in $\End (V)$ by the formula:
\[ \{a, b\} = ab - (-1)^{\tilde a \tilde b}ba
\]
It is easy to check that $\left(\End (V), \{\cdot, \cdot\}\right)$
is a Lie superalgebra. 

\hfill

\remark
Given a $\Z$-graded vector space $A$, one defines
$A^\ev$ as the direct sum of even components, and $A^\odd$ as
the direct sum of odd components. Then a $\Z$-graded
Lie superalgebra is given by a supercommutator on $A$
satisfying $\{A^p, A^q\}\subset A^{p+q}$ and satisfying the
graded Jacobi identity. In the sequel, all  Lie superalgebras
we consider are of this type. An endomorphism $u\in \End(A)$
is called {\bf even} if $u(A^\odd)\subset A^\odd$ and
$u(A^\ev)\subset A^\ev$, and {\bf odd} if
$u(A^\ev)\subset A^\odd$ and
$u(A^\odd)\subset A^\ev$. An endomorphism which is either
odd or even is called {\bf pure}.

\hfill

\definition
A graded algebra $A$ is called {\bf graded commutative}
if $\{a, b\} =0$ for all $a, b\in A$.

\hfill

The Grassmann algebra and de Rham algebra are clearly graded commutative.

\hfill

\definition
Let $\g$ be a graded commutative algebra.
A map $\delta:\; \g\arrow \g$ is called an {\bf even derivation}
if it is even 
and satisfies $\delta(xy)= \delta(x)y+ x\delta(y)$.
It is called an {\bf odd derivation} if it is odd and satisfies
$\delta(xy)= \delta(x)y+ (-1)^{\tilde x}x\delta(y)$.
It is called {\bf graded derivation}, or {\bf superderivation}, if 
it shifts the grading by $i$ and satisfies
$\delta(ab) = \delta(a) b + (-1)^{ij} a \delta(b)$,
for each $a \in A^j$. 

\hfill

\remark
The supercommutator of two superderivations is again a derivation.
Therefore, the derivations form a Lie superalgebra.

\subsection{Differential operators on graded commutative algebras}
\label{_diffe_super_alge_Subsection_}

We need some algebraic results, which are almost trivial,
and well-known for commutative algebras. We extend these statements
to graded commutative algebras; the proofs are the same
as in the commutative setting (see \cite{_co_}).

\hfill

\definition\label{_diff_ope_Definition_}
Let $A^*$ be a graded commutative algebra.
{\bf The algebra of differential operators} $\Diff(A^*)$ is an associative
subalgebra of $\End(A^*)$ generated by
graded derivations and $A^*$-linear self-maps.
Let $\Diff^0(A^*)$ be the space of $A^*$-linear self-maps,
$\Diff^0(A^*)=A^*$, $\Diff^1(A^*)\supset \Diff^0(A^*)$
the subspace generated by $\Diff^0(A^*)$ and all graded derivations, and
$\Diff^i(A^*):= \Diff^{i-1}(A^*)\cdot \Diff^1(A^*)$.
This gives a multiplicative filtration
$\Diff^0(A^*)\subset \Diff^1(A^*) \subset \Diff^2(A^*)\subset \cdots$ 
The elements of $\Diff^i(A^*)$ are called 
{\bf differential operators of order $i$} on $A^*$.

\hfill

\claim\label{_diff_ope_commutator_Claim_}
Let $D\in \Diff^i(A^*), D'\in \Diff^j(A^*)$ be differential
operators on a graded commutative algebra. Then
$\{D, D'\} \in \Diff^{i+j-1}(A^*)$. 

\proof Since the 
commutator of two derivations is a derivation, one has
$\{\Diff^1(A^*), \Diff^1(A^*)\}\subset \Diff^1(A^*)$.
Then we use induction on $i$ and the standard
commutator identities in the associative algebra. 
\endproof

\hfill

We shall apply this claim to geometric operations on
the de Rham algebra, obtaining differential operators of first order.
Note that the ``differential operator'' in the usual sense
is a different notion. For example, the  interior product operator $i_v$  of contraction 
with a vector field $v$ on a manifold $M$ is an odd derivation of the de Rham algebra, 
hence it is a first order differential operator; however,
$i_v$ is $C^\infty(M)$-linear, and thus it is not a differential operator in the usual sense.

\hfill

\claim\label{_first_order_vs_deriva_Claim_}
Let $D$ be a differential operator of first order on $A^*$.
Then $D(x)= D(1) x + \delta(x)$, where $\delta$ is a derivation.

\proof This is seen by defining $\delta:=D-D(1)\cdot$, then observing that $\delta(1)=0$, which implies that $\delta\circ a-a\circ\delta=\delta(a)$ (we identify $a\in A$ with $a\in\End(A)$, $a(b)=ab$). We then have
$$\delta(ab)=(\delta\circ a)(b)=(a\circ\delta)(b)+\delta(a)(b)=a\delta(b)+\delta(a)b,$$
proving that $\delta$ is a derivation.
\endproof

\hfill

We use this formalism to compare first order 
differential operators on $A^*$ as follows.

\hfill

\claim\label{_deriva_dete_Claim_}
Let $\delta$ be a derivation on $A^*$. Then
$\delta$ is uniquely determined by the values it takes
on any set of multiplicative generators of $A^*$.
\endproof

\hfill

\claim\label{_first_or_dete_Claim_}
Let $D$ be a first order differential operator on $A^*$. Then
$D$ is uniquely determined by $D(1)$ and the values it takes
on any set of multiplicative generators of $A^*$

\proof By \ref{_first_order_vs_deriva_Claim_}, 
$D-D(1)$ is a derivation, hence \ref{_deriva_dete_Claim_} implies 
\ref{_first_or_dete_Claim_}.
\endproof

\hfill

\corollary\label{_d^2=0_dete_Corollary_}
Let $d_1, d_2:\; \Lambda^*(M) \arrow \Lambda^{*+1}(M) $ be first order 
differential operators\footnote{Here, as elsewhere, ``differential operators
on the de Rham algebra'' are understood in the algebraic sense as above.} 
on  the de Rham algebra of a manifold $M$,
satisfying $d_1^2=d_2^2=0$, and $V\subset \Lambda^1(M)$ a subspace such that
the space $d_1(C^\infty (M)) + V$  generates $\Lambda^1(M)$
as a $C^\infty (M)$-module. Suppose that $d_1\restrict {C^\infty (M)}=d_2\restrict {C^\infty (M)}$
and $d_1\restrict {V}=d_2\restrict {V}$. Then $d_1=d_2$.

\proof Clearly, $d_1=d_2$ on $C^\infty (M)\oplus d_1(C^\infty (M))$.
Then $d_1=d_2$ on the set of multiplicative generators 
$C^\infty (M)+ d_1(C^\infty (M))+ V$, and \ref{_first_or_dete_Claim_}
implies \ref{_d^2=0_dete_Corollary_}. \endproof

\hfill

The following claim is used many times in the sequel.

\hfill

\claim\label{_d^2=0_super_Claim_}
Let $d$ be an odd element in a Lie superalgebra $\goth h$,
satisfying $\{d,d\}=0$. Then
$\{d, \{d, u\}\}=0$ for any $u\in {\goth h}$.

\proof By the super Jacobi identity,
\begin{equation}\label{_odd_squa_Jaco_Equation_}
\{d, \{d, u\}\} = - \{d, \{d, u\}\} + \{\{ d, d\}, u\} = - \{d, \{d, u\}\}.
\end{equation}
\endproof

This claim is a special case of the following: 

\hfill

\claim\label{_d^2neq0_super_Claim_}
Let $d$ be an odd element in a Lie superalgebra $\goth h$. Then
$2\{d, \{d, u\}\}= \{\{d,d\}, u\}$ for any $u\in {\goth h}$.

\proof Follows from \eqref{_odd_squa_Jaco_Equation_}.
\endproof

\subsection{Supersymmetry on K\"ahler manifolds}
\label{_susy_Kah_Subsection_}

One of the purposes of this paper is to obtain a natural superalgebra
acting on the de Rham algebra of a Sasakian manifold.
This is modeled on the superalgebra of a K\"ahler manifold,
generated by the de Rham differential, Lefschetz triple,
and other geometric operators. To make the analogy
more clear, we recall the main results on the supersymmetry 
algebra of K\"ahler manifolds. We follow \cite[Section 1.3]{_Verbitsky:HKT_}.

\hfill

Let $(M, I, g, \omega)$ be a K\"ahler manifold. Consider $\Lambda^*(M)$ as a graded
vector space. The differentials $d, d^c:= - I dI=IdI^{-1}$
can be interpreted as odd elements in $\End(\Lambda^*(M))$,
and the Hodge operators $L, \Lambda, H$ as even elements. 
As usual, we denote the supercommutator 
as $\{\cdot,\cdot\}$. In terms of the associative
algebra, $\{a,b\} = ab +ba$ when $a, b$ are odd,
and $\{a,b\} = ab -ba$ if at least one of them is even.
Let $d^*:=  \{\Lambda, d^c\}$, $(d^c)^*:= -\{\Lambda, d\}$.
The usual Kodaira relations can be stated as follows
\begin{equation}\label{_Kodaira_Equations_}
\begin{aligned} 
 {} & \{ L, d^*\} = - d^c,\ \ \ \  \{ L, (d^c)^*\} = d, \ \ \ \ 
\{ d, (d^c)^*\} = \{ d^*, d^c\} =0, \\
 {} & \{ d, d^c\} =  \{ d^*, (d^c)^*\} =0,\ \ \ \  
\{ d, d^*\} = \{ d^c, (d^c)^*\} = \Delta,
\end{aligned}
\end{equation}
where $\Delta$ is the Laplace operator,
 commuting with $L, \Lambda, H$, and
$d$, $d^c$.

\hfill

\definition \label{_KdR_Definition_}
Let $(M, I, g, \omega)$ be a K\"ahler manifold.
Consider the Lie superalgebra $\a \subset \End(\Lambda^*(M))$
generated by the following operators:
\begin{enumerate}
\item $d$, $d^*$, $\Delta$, constructed out of the Riemannian metric.

\item $L(\alpha):= \omega\wedge \alpha$.

\item $\Lambda(\alpha) := * L * \alpha$.  It is easily seen that $\Lambda= L^*$.

\item The Weil operator $W\restrict{\Lambda^{p,q}(M)}=\1(p-q)$
\end{enumerate}
This Lie superalgebra is called {\bf the algebra of supersymmetry
of the K\"ahler manifold}.

\hfill

Using \ref{_kah_susy_Theorem_} below,
it is easy to see that
$\goth a$ is in fact independent from $M$.

\hfill

This Lie superalgebra was studied from the physicists'
point of view in \cite{_FKS_}.

\hfill

\theorem\label{_kah_susy_Theorem_}
Let $M$ be a K\"ahler manifold, and 
$\goth a$ its  supersymmetry algebra
acting on $\Lambda^*(M)$. Then
$\goth a$ has dimension $(5|4)$
(that is, its odd part is 4-dimensional, and its even part
is 5-dimensional). The odd part is generated by
$d, d^c=IdI^{-1}, d^*, (d^c)^*$, the even part
is generated by the Lefschetz triple
$L, \Lambda, H=[L, \Lambda]$, the Weil operator $W$ and the
Laplacian $\Delta=\{d, d^*\}$.
Moreover, the Laplacian $\Delta$ is
central in $\goth a$, hence $\goth a$ also acts on the
cohomology of $M$. The following are the only non-zero
commutator relations in $\goth a$:
\begin{enumerate}
\item ${\goth {sl}}(2)$-relations in $\langle L, \Lambda, H=[L, \Lambda]\rangle$:
\[
[H, L]=2L, \ \  [H, \Lambda]=-2\Lambda,\ \  [L, \Lambda]=H.
\]
For any operator $D$ of grading $k$, one has $[H, D]=kD$.

\item The Weil operator acts as a complex structure on the odd part of $\goth a$:
\[ [W, d]=d^c, \ \   [W, d^c]=-d, \ \ [W, d^*]=-(d^c)^*, \ \   [W, (d^c)^*]=d^*.\]
\item The K\"ahler-Kodaira relations between the differentials and the Lefschetz operators are:
\begin{equation}\label{_Kodaira_rel_Equation_}
  [\Lambda, d] = (d^c)^*,\ \  [ L, d^*] = - d^c,\ \ [\Lambda, d^c] = - d^*, \ \ [ L, (d^c)^*] = d.
\end{equation}
\item Almost all odd elements supercommute, with the only exception
\[
\Delta=\{d, d^*\}=\{d^c, (d^c)^*\},
\]
and $\Delta$ is central. In other words, the odd elements of $\goth a$
generate the odd Heisenberg superalgebra, see \ref{_odd_Heisenberg_Claim_}.
\end{enumerate}

\proof
These relations are standard in algebraic geometry (see e.g. \cite{_Griffi_Harri_}), but probably the easiest way to prove them is using the results about Lie superalgebras
collected in Subsection \ref{_diffe_super_alge_Subsection_}.

\hfill

{\bf Proof of the Lefschetz ${\goth {sl}}(2)$-relations:}
These relations would follow if we prove that
$H:=[L, \Lambda]$ acts on $p$-forms by multiplication
by $p-n$, where $n=\dim_\C M$. Since $L, \Lambda, H$ are
$C^\infty(M)$-linear, it would suffice to prove these
relations on a Hermitian vector space. 

Let $V$ be a real vector space equipped with
a scalar product, and fix an orthonormal basis $\{v_1, \ldots, v_{m}\}$.
Denote by $e_{v_i}:\; \Lambda^k V  \arrow \Lambda^{k+1} V$
the operator of multiplication, $e_{v_i}(\eta) = v_i \wedge \eta$.
Let $i_{v_i}:\; \Lambda^k V  \arrow \Lambda^{k-1} V$
be the operator of contraction with $v_i$.
The following claim is clear.

\hfill

\claim\label{_odd_Heisenberg_Claim_}
The operators $e_{v_i}$, $i_{v_i}$, $\Id$ form a basis of the {\bf 
odd Heisenberg Lie superalgebra}, 
with the only non-trivial 
supercommutator given by the formula $\{ e_{v_i}, i_{v_j}\} = \delta_{ij}\Id$.
\endproof

\hfill

Let now $V$ be an even-dimensional real vector space equipped with
a scalar product, and $\{x_1, \ldots, x_n, y_1, \ldots, y_{n}\}$ an orthonormal basis.
Consider the complex structure operator $I$ such that
$I(x_i)=y_i$, $I(y_i)=-x_i$. 
The fundamental symplectic form is given
by $\sum_i x_i \wedge y_i$, hence
\[
L= \sum_i e_{x_i} e_{y_i}, \ \ \ \Lambda= \sum_i i_{x_i} i_{y_i}.
\]
Clearly, for any odd elements $a, b, c, d$ such that
$\{a,b\}= \{a, d\} = \{b, c\}=\{c,d \}=0$, 
one has
$\{ab, cd\}= -\{a, c\}bd + ca \{b,d\}$. 
Then 
\begin{equation*}
\begin{split}
 [L, \Lambda] &= 
\left [\sum_i e_{x_i} e_{y_i}, \sum i_{x_i} i_{y_i}\right]
=  \sum_{i=1}^{n} e_{y_i} i_{y_i} -\sum_{i=1}^{n}i_{x_i} e_{x_i}\\
& =
\sum_{i=1}^{n} e_{y_i} i_{y_i} +\sum_{i=1}^{n}(e_{x_i} i_{x_i} -1)
\end{split}
\end{equation*}
This term, applied to a monomial $\alpha$ of degree $d$,
would give $(d-n)\alpha$. This proves the Lefschetz ${\goth {sl}}(2)$-relations.

\hfill

{\bf Proof of the relations between the Weil operator $W$ and the odd part of $\goth a$.}
Clearly, it is enough to prove $[W, d]=d^c$, the remaining relations follow by duality or
by complex conjugation. Writing the Hodge components of $d= d^{1,0}+d^{0,1}$, with
$d^{1,0}=\frac{d + \1 d^c}2$ and $d^{0,1}=\frac{d -\1 d^c}2$ ,
we obtain $[W, d]=\1 d^{1,0}-\1 d^{0,1}= d^c$. 

\hfill

{\bf Proof of the K\"ahler-Kodaira relations between the Lefschetz ${\goth {sl}}(2)$-operators and the
odd part of $\goth a$.} As before, it is enough to prove $[ L, d^*] = - d^c$,
the remaining K\"ahler-Kodaira relations follow by duality or
by complex conjugation. The operator $L$ is $C^{\infty}(M)$-linear, hence
it is a differential operator of order 0. The operator
$d^*$ can be written in a frame $\{v_i \}$ of $TM$ 
as 
\begin{equation} \label{_d^*_nabla_Equation_}
d^*(\eta)= \sum_i i_{v_i}\nabla_{v_i}\eta, 
\end{equation}
where $\nabla$ is the Levi-Civita connection of the metric $g$.
Since $\nabla_{v_i}$ and $i_{v_i}$ are
both derivations of the de Rham algebra, their product is
an order 2 differential operator (in the algebraic sense
as given by \ref{_diff_ope_Definition_}). From 
\ref{_diff_ope_commutator_Claim_} it follows that
$[ L, d^*]$ is a first order operator.

We shall prove $[ L, d^*] = - d^c$ by applying 
 \ref{_d^2=0_dete_Corollary_}. First, let us show that
$[L,[L, d^*]]=0$. Clearly, $[\Lambda, d^*]= ([L, d])^*=0$,
and $[H, d^*]=-d^*$, hence $d^*$ is the lowest weight vector in
a weight 1 representation of $\goth{sl}(2)$. This gives
$[L,[L, d^*]]=0$. Now, the super Jacobi identity gives
\[
\{\{L, d^*\},\{L, d^*\}\}= \{L,\{d^*,\{L, d^*\}\}+ \{d^*,\{L,\{L, d^*\}\}\}\}.
\]
The first term in the RHS vanishes by \ref{_d^2=0_super_Claim_}, and 
the second term vanishes because $[L,[L, d^*]]=0$.
Then $([L, d^*])^2=0$. Clearly, $d^c(C^\infty (M))$ generates $\Lambda^1(M)$ over $C^\infty (M)$.
To deduce $[ L, d^*] = - d^c$ from  \ref{_d^2=0_dete_Corollary_}, it remains to show
 that $[ L, d^*]\restrict{C^\infty (M)} = - d^c\restrict{C^\infty (M)}$.
This is clear for the following reason. For any function
$f\in C^\infty (M)$, one has $[ L, d^*] f= - d^*(f \omega)$. Writing
$d^*= \sum_i i_{v_i}\nabla_{v_i}$ as in \eqref{_d^*_nabla_Equation_},
and using $\nabla\omega=0$, we obtain
\[ 
d^*(f \omega)= \sum_i i_{v_i}\nabla_{v_i}(f\omega) = \sum_i \Lie_{v_i}(f) i_{v_i}(\omega)=
 \sum_i - \Lie_{v_i}(f) I(v_i)^\flat= -d^cf.
\]
This finishes the proof of $[ L, d^*] = - d^c$.

\hfill

{\bf Proof of the commutator relations between the odd part of $\goth a$.}
We have already shown that $d^c=[W, d]$. Then $\{d, d^c\}= \{d,\{W, d\}\}=0$
by \ref{_d^2=0_super_Claim_}. Similarly, $d^*= - [\Lambda, d^c]$, giving
$\{d^c, d^*\}=0$. The relation $\{(d^c)^*, d\}=0$ is obtained by duality.
Finally, $\{d, d^*\}=\{d^c, (d^c)^*\}$ is obtained by
applying $[\Lambda, \cdot]$ to $\{d, d^c\}=0$. Using 
the K\"ahler-Kodaira relations \eqref{_Kodaira_rel_Equation_}, we obtain
\[
0=\{\Lambda, \{d, d^c\}\}= \{\{\Lambda,d\}, d^c\}+ \{d, \{\Lambda,d^c\}\}=
\{(d^c)^*,  d^c\}-\{d, d^*\},
\]
giving $\{d, d^*\}=\{d^c, (d^c)^*\}$.

We proved all the relations in the K\"ahler supersymmetry algebra
$\goth a$, finishing the proof of \ref{_kah_susy_Theorem_}.
\endproof

\hfill

Further in this paper, we shall develop similar relations
for the superalgebra associated with a Sasakian manifold.


\section{Sasakian manifolds: definition and the basic
  notions}
\label{_Sasak_Section_}

In this section we provide the necessary background on Sasakian manifolds. For details, please see \cite{_Blair_}, \cite{_Boyer_Galicki_}. The most convenient definition for our context is the  one which relates Sasakian and  K\"ahler geometries in terms of Riemannian cones.

\hfill

\definition  A {\bf Sasakian manifold} is a Riemannian manifold $(S, g)$ with
	a K\"ahler structure on its Riemannian cone $C(S):=(S\times \R^{>0}, t^2g+dt^2)$, such that the
	homothety map $h_\lambda:\; C(S)\arrow C(S)$ mapping
	$(m,t)$ to $(m,\lambda t)$ is holomorphic.

\hfill

\remark \label{_structure_of_cone_rem_}
(i) A Sasakian manifold is clearly contact, because its cone
is symplectic and $h_\lambda$ acts by symplectic homotheties.

(ii) Let $S$ be a Sasakian manifold, $\omega$  the K\"ahler form
on $C(S)$, and  $\xi=t\frac{d}{dt}$ the homothety vector
field along the generators of the cone (it is also called
{\bf Euler field}). The contact form is $\eta=i_\xi\omega\restrict {t=1}$. 
Then $\Lie_{I\xi}t= \langle dt, I\xi\rangle=0$ (where $\Lie_v$ denotes the Lie derivative in the direction of $v$), and  hence $I\xi$
is tangent to $S\subset C(S)$. Then 
$\eta(I\xi)=1$ and $i_{I\xi}  d\eta=0$, and thus {\bf the
	Reeb vector field of a Sasakian manifold is $\vecr=I\xi$}. 
 
 (iii) The function 
 $t^2$ is a K\"ahler potential on $C(S)$. Moreover,
 the form $dd^c \log t$
 vanishes on $\langle \xi, I(\xi)\rangle$ and the rest of its
 eigenvalues are positive.

\hfill

\proposition \label{_Reeb_Sasakian_Theorem_}
The Reeb field has constant length and acts on a Sasakian manifold
by contact isometries, i.e.  the flow of
$\vecr$ contains only isometries which preserve the
contact subbundle. Moreover, its action lifts
to a holomorphic action on its cone.

\hfill

\remark\label{_Volume_form_on__Sasaki_Theorem_}
On a $2n+1$-dimensional Sasakian manifold, the contact form satisfies $\eta\wedge(d\eta)^n\neq 0$. Since $i_{\vecr}d\eta=0$, we see that $d\eta$ is a volume form on the contact distribution $\vecr^\perp$.

\hfill 

\definition A Sasakian (resp. contact)  manifold is called {\bf regular} if its 
Reeb field generates a free action of $S^1$; it is called {\bf 
	quasi-regular} if all orbits of the Reeb field are closed, and it is called 
{\bf irregular} otherwise.

\hfill

\example \label{_ample_line_bdle_ex_} Let $S$ be a  regular Sasakian manifold, and $\vecr$ its Reeb
field. Then the space of $\vecr$-orbits $X$ is K\"ahler.
Moreover, $X$ is equipped with a positive holomorphic
Hermitian  line bundle
$L$ such that $S$ is the space of unit vectors in $L$. Conversely: if $X$ is a compact projective manifold, together with an ample line bundle $L\arrow X$, then the space of unit vectors in $L$ is a regular Sasakian manifold.


\section{Hattori differentials on Sasakian manifolds}
\label{_Hattori_Section_}


\subsection{Hattori spectral sequence and associated differentials}
\label{_Leray_Serre_Subsection_}

Let $\pi:\; M \arrow B$ be a smooth fibration,
and $F_k\subset \Lambda^*(M)$ be the  ideal generated
by $\pi^*\Lambda^k(B)$. The Hattori spectral sequence
(\cite{_Hattori_}) is the spectral sequence 
associated with this filtration.

The $E_1^{p,q}$-term of this sequence is $\Lambda^p(B)\otimes_\R R^q\pi_* \R_M$,
where $R^*\pi_* \R_M$ denotes the local system of cohomology of the fibres, and
the $E_2^{p,q}$-term is $H^p(B, R^q\pi_* \R_M)$.

The same can be done when $M$ is a manifold equipped
with an integrable distribution $F\subset TM$, giving a
spectral sequence converging to $H^*(M)$.
This is done as follows.

\hfill

\definition
Let $M$ be a manifold, and $F\subset TM$ an
integrable distribution. A $k$-form $\alpha\in \Lambda^* (M)$
is called {\bf basic} if for any vector field $v\in F$,
one has $\Lie_v\alpha=0$ and $i_v\alpha=0$.

\hfill

If $F$ is the tangent bundle of the fibres of a
fibration $\pi:\; M \arrow B$, then the space
of basic forms is $\pi^*\Lambda^k(B)$.
We are going to produce a spectral sequence which gives
the standard Hattori spectral sequence when $F$ 
is tangent to the leaves of a fibration.

Define {\bf the Hattori filtration associated with $F$} as   
$F_n\subset F_{n-1}\subset \cdots$ 
by putting $F_k\subset \Lambda^*(M)$, where $F_k$ is the
ideal generated by basic $k$-forms.

When $M$ is Riemannian and the metric is compatible with the foliation, 
this filtered bundle is decomposed 
into the direct sum of subquotients. This
gives a decomposition of the de Rham differential,
$d= d_0 + d_1 + d_2+\cdots + d_{r+1}$, where $r=\rk F$
with each successive piece associated with the differential
in $E_k^{p,q}$, as follows.

Using the metric, we split the cotangent bundle into orthogonal complements
as $\Lambda^1(M)= \Lambda^1_\hor(M) \oplus \Lambda^1_\vert(M)$,
where $\Lambda^1_\hor(M)$ is generated by basic 1-forms, and
$\Lambda^1_\vert(M)= F^*$ is its orthogonal complement.
Denote by $\Lambda^p_\hor(M)$, $\Lambda^q_\vert(M)$
the exterior powers of these bundles.

This gives the following splitting of the de Rham algebra of 
$M$:
\begin{equation}\label{_horis_vert_decompo_Equation_}
\Lambda^m(M)= \bigoplus_p \Lambda^p_\hor(M)\otimes \Lambda^{m-p}_\vert(M) 
\end{equation}
with $F_{p}/F_{p-1}= \bigoplus_q \Lambda^p_\hor(M)\otimes \Lambda^{q}_\vert(M)$
where $F_*$ denotes the Hattori filtration.

Consider the associated
decomposition of the de Rham differential,
$d= d_0 + d_1 + d_2+ \cdots + d_{r+1}$, where $r=\rk F$,
and 
\[ d_i:\; \Lambda^p_\hor(M)\otimes \Lambda^{q}_\vert(M)\arrow 
   \Lambda^{p+i}_\hor(M)\otimes \Lambda^{q+1-i}_\vert(M)
\]
The terms $d_{i}$ vanish for $i> r+1$ because
$F$ is $r$-dimensional, hence for $i>r+1$ either $\Lambda^{q+1-i}_\vert(M)$
or $\Lambda^{q}_\vert(M)$ is 0.

These differentials are related to the differentials in
the Hattori spectral sequence in the following way:
to find $E_1^{p,q}$, one takes the cohomology of $d_0$.
Then one restricts $d_1$ to $E_1^{p,q}$, and its cohomology
gives $E_2^{p,q}$, and so on. In this sense, the
Hattori differentials are indeed differentials
in the Hattori spectral sequence.

\hfill

\remark 
Each of the differentials $d_i$ is a derivation, because
the decomposition $\Lambda^*(M)= \bigoplus_{p,q}\Lambda^p_\hor(M)\otimes \Lambda^q_\vert(M)$
is multiplicative.

\hfill

The spectral sequence which we call ``Hattori spectral sequence''
was re-invented independently on several instances after
Hattori. In a book \cite[Section 1.6]{_Brylinski_} by J.-L. Brylinski
it was described as ``Cartan spectral sequence'', without 
reference. About the same time, it was described
in Vlad Sergiescu's Ph. D. thesis  and in his subsequent papers
(\cite{_Sergiescu:thesis_,_El_Kacimi_Sergiescu_Hector_,_Sergiescu:Mexican_})
under the name ``Leray-Serre type spectral sequence''.
This work was quite influential, with a number of publications
citing Sergiescu's papers and his thesis 
(for example, \cite{_Alvares:finiteness_,_Alvares:duality_,_Domingues:finiteness_}). 
In \cite{_Alvares:finiteness_},
it was shown that all terms on $E_2$ page of this spectral
sequence are finite-dimensional when the foliation
admits a transversal Riemannian structure, and
in \cite{_Domingues:finiteness_} the same result
was proven for cohomology with coefficients in a local system.

\subsection{Hattori differentials on Sasakian
  manifolds}
\label{_LS_decompo_Sasakian_Subsection_}

Let now $Q$ be a Sasakian manifold, $\vecr$ its Reeb field, normalized
in such a way that $|\vecr|=1$, and
$R\subset TQ$ the 1-dimensional foliation generated by the Reeb field.
The corresponding Hattori differentials are written as
$d=d_0 + d_1 + d_2$, because $R$ is 1-dimensional. 
Since $d^2=0$, one has $d_0^2=d_2^2=0$ and
$\{d_0, d_2\}=-\{d_1, d_1\}$. The differentials $d_0, d_2$
can be described explicitly as follows.

\hfill

\claim\label{_d_0_expli_Claim_}
Let $e_\vecr:\; \Lambda^*(Q)\arrow \Lambda^{*+1}(Q)$ be the operator of multiplication
by the form $\vecr^{\,\flat}=\eta$ dual to $\vecr$. Then $d_0= e_\vecr \Lie_\vecr$.

\proof 
Locally, the sheaf $\Lambda^1(Q)$ is generated over $C^\infty (Q)$ by 
the basic 1-forms and $\vecr^{\,\flat}$. Clearly, the differential of a basic
1-form $\alpha$ belongs to $\Lambda^2_\bas(Q)$, hence $d_0(\alpha)=0$,
and $e_\vecr \Lie_\vecr(\alpha)=0$. Also, 
$d_0(\vecr^{\,\flat})= e_\vecr \Lie_\vecr(\vecr^{\,\flat})=0$. 
By \ref{_deriva_dete_Claim_}, to prove $d_0=e_\vecr \Lie_\vecr$
it remains to show that these two operators are equal on 
$C^\infty (Q)$. However, on $C^\infty (Q)$, 
we have $d_0= e_\vecr \Lie_\vecr$, because $d_0f$ is the orthogonal
projection of $df$ to $\Lambda^1_\vert(Q)=R^*$ generated
by $e_\vecr$, for all $f\in C^{\infty}(Q)$.
\endproof

\hfill

Recall that the space of leaves of $R$ on a regular Sasakian manifold is equipped with 
a complex and K\"ahler structure, \ref{_ample_line_bdle_ex_}. The corresponding 
K\"ahler structure can be described very explicitly.

\hfill

\claim\label{_omega_0_trans_Kahler_Claim_}
Let $Q$ be a regular Sasakian manifold, and $\vecr^{\,\flat}$ its contact form.
Then $\omega_0:=d(\vecr^{\,\flat})$ is basic with respect to $R$,
and defines a transversally K\"ahler structure, that is, the 
K\"ahler structure on the space of leaves of $R$ (\ref{_transvers_Kahler_Definition_}).

\hfill

\proof Let $X=Q/\vecr$ be the space of orbits. This quotient is well
defined and smooth, because $\vecr$ is regular.
Then $X=C(Q)/\C^*$, where the $\C^*$-action is generated
by $ \xi=t\frac{d}{dt}$ and $\vecr= I(\xi)$, and hence it
is holomorphic. 
Therefore,
	$X$ is a complex manifold (as  a quotient of a complex
manifold by a holomorphic action of a Lie group).
It is K\"ahler by \ref{_structure_of_cone_rem_} (iii). \endproof

\hfill

\proposition\label{_d_2_expli_Proposition_}
Let $L_{\omega_0}:\; \Lambda^*(Q)\arrow \Lambda^{*+2}(Q)$ be the operator of multiplication
by the transversally K\"ahler form $\omega_0=d(\vecr^{\,\flat})$, and
$i_\vecr$ the contraction with the Reeb field. Then
$d_2 = L_{\omega_0} i_\vecr$.

\hfill

\proof Clearly, the 
Hattori differentials 
$$d_i:\; \Lambda^p_\hor(Q)\times \Lambda^{q}_\vert(Q)\arrow 
   \Lambda^{p+i}_\hor(Q)\otimes \Lambda^{q+1-i}_\vert(Q)$$
vanish on $\Lambda^0(Q)$ unless $i=0,1$. Therefore, the differentials
$d_2, d_3, ...$ are always $C^\infty(Q)$-linear.
By \ref{_deriva_dete_Claim_}, it only remains to show  that $d_2 = L_{\omega_0} i_\vecr$ on 
some set of 1-forms generating $\Lambda^1(Q)$ over $C^\infty (Q)$.

Clearly, on $\Lambda^1_\hor(Q)$ the differential
$d_2$  should act as 
$$d_2:\;\Lambda^1_\hor(Q)\arrow \Lambda^{3}_\hor(Q)\otimes \Lambda^{-1}_\vert(Q),$$
hence $d_2\restrict{\Lambda^1_\hor(Q)}=0$.
To prove \ref{_d_2_expli_Proposition_} it remains to show that 
$d_2(\vecr^{\,\flat})=L_{\omega_0} i_\vecr(\vecr^{\,\flat})$.
However, $d_2(\vecr^{\,\flat})$ is the $\Lambda^2_\hor(Q)$-part of  $d(\vecr^{\,\flat})$,
giving $d_2(\vecr^{\,\flat})=d(\vecr^{\,\flat})=\omega_0$, and $L_{\omega_0} i_\vecr(\vecr^{\,\flat})=\omega_0$
because $i_\vecr(\vecr^{\,\flat})=1$.
\endproof

\hfill

The Hattori differential $d_1$ is, heuristically speaking, the ``transversal component''
of the de Rham differential. Indeed, $d_1(\eta)=d(\eta)$ for any basic form $\eta$.
Since the leaf space of $R$ is equipped with a complex structure, it is natural to expect that
the Hodge components of $d_1$ have the same properties as the Hodge components
of the de Rham differential on a complex manifold.

\hfill

\claim\label{_Hodge_decompo_d_1_Sasa_Claim_}
Let $Q$ be a Sasakian manifold, and 
$\Lambda^m(Q)= \bigoplus_p \Lambda^p_\hor(Q)\otimes \Lambda^{m-p}_\vert(Q)$ the decomposition
associated with $R\subset TQ$ as in \eqref{_horis_vert_decompo_Equation_}.
Using the complex structure on the basic forms, consider the
Hodge decomposition $\Lambda^m_\hor(Q)=\bigoplus_p \Lambda^{m-p,p}_\hor(Q)$.
Then the differential $d_1:\; \Lambda^p_\hor(Q)\otimes \Lambda^{q}_\vert(Q)\arrow 
\Lambda^{p+1}_\hor(Q)\otimes \Lambda^{q}_\vert(Q)$ 
has two Hodge components $d_1^{0,1}$ and $d_1^{1,0}$.
Moreover, the differential $d_1^c:=d_1^{0,1}-d_1^{1,0}$
satisfies $d_1^c=I d_1 I^{-1}$, where $I$ acts as
$\1^{p-q}$ on $\Lambda^{p,q}_\hor(Q)\otimes \Lambda^m_\vert(Q)$.

\proof
First, let us prove that $d_1$ has only two non-zero Hodge components.
A priori, $d_1$ could have several Hodge components,
$d_1= d^{-k, k+1}_1+ d^{-k+1, k}_1+ \cdots + d^{k, -k+1}_1+d^{-k,
  k+1}_1$. 
This is what happens with the Hodge components of the  de Rham
differential on an almost complex manifold.
All these components are clearly derivations.
However, 
\[ 
d(\Lambda^0(Q)) \subset \Lambda^{1}_\vert(Q)\oplus
\Lambda^{1,0}_\hor(Q) \oplus\Lambda^{0,1}_\hor(Q). 
\]
Therefore, only $d_1^{0,1}$ and $d_1^{1,0}$ are non-zero on
functions, the rest of the Hodge components are $C^\infty(Q)$-linear.
By \ref{_deriva_dete_Claim_}, it would suffice to show that
the other differentials vanish on a basis in $\Lambda^1(Q)$.

Since the space of leaves of $R$ is a complex manifold,
and $d=d_1$ on basic forms, we have $d_1=d_1^{0,1}+d_1^{1,0}$
on basic forms. Since $d(\vecr^{\,\flat})=\omega_0$, we find 
$d_1(\vecr^{\,\flat})=0$, which gives $d_1=d_1^{0,1}+d_1^{1,0}$
on $\vecr^{\,\flat}$. We proved the decomposition
$d_1=d_1^{0,1}+d_1^{1,0}$. The relation
$d_1^{0,1}-d_1^{1,0}=I d_1 I^{-1}$ follows in the usual way,
because $\frac{d_1+ \1 d_1^c}{2}$ has Hodge type $(1,0)$,
hence satisfies $d_1^{1,0}=\frac{d_1+ \1 d_1^c}{2}$.
\endproof


\section{Transversally K\"ahler manifolds}
\label{_Tra_Kah_Section_}


\definition\label{_transvers_Kahler_Definition_}
A manifold $M$ equipped with an integrable distribution
$F\subset TM$ is called a {\bf foliated manifold}. In the sequel, we shall
always assume that $F$ is orientable. 
Let $\omega_0\in \Lambda^2 (M)$ be a closed, basic 2-form 
on a foliated manifold $(M,F)$, vanishing on $F$ and
non-degenerate on $TM/F$. Let $g_0 \in \Sym^2(T^*M)$ 
be a basic bilinear symmetric form which is positive definite on $TM/F$.
Since $g_0, \omega$ are basic, the operator $I:= \omega^{-1}_0 \circ g_0 :\; TM/F\arrow TM/F$
is well defined on the leaf  space $L$ of $F$ (locally the leaf space always
exists by Frobenius theorem).  Assume that $I$ defines an integrable complex structure
on $L$ for any open set $U\subset M$ for which the leaf space is well defined.
Then $(M, F, g, \omega)$ is called {\bf transversally K\"ahler}.
A vector field $v$ such that $\Lie_v(F)\subset F$, and
$\Lie_v I =0$ is called {\bf transversally holomorphic};
it is called {\bf transversally Killing} if, in addition,
$\Lie_v g_0 =0$.

\hfill

\definition
Let $F\subset TM$ be an integrable distribution and $\Lambda_\bas^*(M)$
the complex of basic forms. Its cohomology algebra is called the 
{\bf basic}, or {\bf transversal} cohomology of $M$. We denote the
basic cohomology by $H_\bas^*(M)$.

\hfill

\remark
Note that $H_\bas^*(M)$ can be infinite-dimensional even
when $M$ is compact, \cite{_Schwarz:foliation_}.

\hfill

The main result of this section is the following theorem,
which is a weaker form of the main theorem from
\cite{_Kacimi_}. Our result is less general,
but the proof is simple and self-contained.

\hfill

\theorem\label{_Transversal_Lefschetz_Theorem_}
Let $(M, F, \omega_0, g_0)$ be a compact, transverally K\"ahler manifold.
Assume, moreover that:
\begin{itemize}
	\item[(*)] $M$ is equipped with a Riemannian metric 
$g$ such that the restriction of $g$ to the orthogonal complement
$F^\bot = TM/F$ coincides with $g_0$, and $F$ is generated by a
collection of Killing vector fields $v_1, ..., v_r$. 
\item[(**)] There exists a closed differential form $\Phi$
  on $M$
 which vanishes on $F$ and gives a Riemannian volume
 form on $TM/F$.\footnote{The assumption (*) and (**) hold
for Sasakian manifolds (\ref{_Reeb_Sasakian_Theorem_} and \ref{_Volume_form_on__Sasaki_Theorem_}) and Vaisman manifolds (\ref{_canonical_foliation_properties_}).}
\end{itemize}
Then the basic cohomology $H^*_\bas(M)$ of $M$ is finite-dimensional and admits
the Hodge decomposition and the Lefschetz ${\goth {sl}}(2)$-action,
 as in the K\"ahler case.

\hfill

\proof Consider the differential graded algebra of basic forms
$\Lambda^*_\bas(M)$.
This algebra is equipped with an action of the superalgebra of
K\"ahler supersymmetry $\goth a$ as in \ref{_kah_susy_Theorem_}. 
Indeed, define the Lefschetz ${\goth {sl}}(2)$-action by taking
the Lefschetz triple $L_{\omega_0}, \Lambda_{\omega_0}:=*L_{\omega_0}*$
and $H_{\omega_0}:= [L_{\omega_0}, \Lambda_{\omega_0}]$, and
the transversal Weil operator $W$ acting in the standard
way on $\Lambda^*_\bas(M)$ and extended to $\Lambda^*(M)$
by acting trivially on $\Lambda^1_\bas(M)^\bot$
(or in any other way, it does not matter). 
Together with the de Rham differential
$d:\; \Lambda^*_\bas(M)\arrow \Lambda^*_\bas(M)$
these operators generate the Lie superalgebra
$\goth a\subset \End(\Lambda^*_\bas(M))$, which
is isomorphic to the superalgebra of
K\"ahler supersymmetry $\goth a$ (\ref{_kah_susy_Theorem_}), 
because $\Lambda^*_\bas(M)$ is locally identified
with the algebra of differential forms on the leaf space of $F$
which is K\"ahler. 

Then \ref{_Transversal_Lefschetz_Theorem_} would follow if
we identify $\ker (\Delta_\bas\restrict {\Lambda^*_\bas(M)})$ with the space 
$H^*_\bas(M)$, where $\Delta_\bas\in \goth a$ is the transversal
Laplace operator, $\Delta_\bas=\{d, d^*_\bas\}$, where
$d^*_\bas$ denotes the $d^*$-operator on the leaf space.

We reduced  \ref{_Transversal_Lefschetz_Theorem_} to the
following result.

\hfill

\proposition\label{_basic_harmo_Proposition_}
Let $(M, F, g)$, $\rk F=r$, be a Riemannian foliated manifold,
$d=d_0+d_1+ \cdots+ d_{r+1}$ the Hattori decomposition of the differential,
and $\Delta_\bas:= \{d, d_\bas^*\}$ the basic Laplacian defined
on basic forms. Assume that:
\begin{itemize}
	\item[(*)] $F$ is generated by a
collection of Killing 
vector fields $v_1, v_2, \ldots$
	\item[(**)] The Riemannian volume form $\Phi\in\Lambda^r_\vert(M)$
satisfies $d_1(\Phi)=0$. 
\end{itemize}
Then there exists a natural
isomorphism between the basic harmonic forms
and $H^*_\bas(M)$.

\hfill

We start from the following lemma.

\hfill

\lemma
In the assumptions of \ref{_basic_harmo_Proposition_}, let
$\alpha, \beta \in \Lambda^*_\bas(M)$ be two basic forms.
Then $g(d\alpha, \beta)= g(\alpha, d^*_\bas\beta)$,
where $d^*_\bas$ denotes the $d^*$-operator on the leaf space of $F$.

\hfill

\proof
This is where we use the assumption (**) of \ref{_basic_harmo_Proposition_}.
Let $d^*_h$ be the composition of $d^*$ with the orthogonal
projection to the horizontal part $\Lambda^*_\hor(M)$
(see Subsection \ref{_Leray_Serre_Subsection_}).
We only need  to show that 
\begin{equation}\label{_d^*_bas_Equation_}
g(d^*_h\alpha, \beta) = g(d^*_\bas\alpha, \beta)\ \ \text{for all}\ \ \alpha, \beta\in \Lambda^*_\bas(M). 
\end{equation}
By \ref{_kah_susy_Theorem_}, one has $d^*_\bas= - *_\bas d *_\bas$
where $*_\bas$ is the Hodge star operator on the leaf space.
Let $r=\rk F$ and $\Phi\in \Lambda^r F$ be the Riemannian volume form.
Using the assumption (**), we obtain that 
$d\Phi\in \bigoplus_{i=0}^{r-1} \Lambda^*_\hor(M) \otimes \Lambda^i_\vert(M)$.
Then $*\alpha= \Phi \wedge *_\bas(\alpha)$, which gives
\[
*d*\alpha = *d( \Phi \wedge *_\bas(\alpha)) = d^*_\bas\alpha + *(d\Phi\wedge \alpha).
\]
The last term belongs to $\bigoplus_{i=1}^{r} \Lambda^*_\hor(M) \otimes \Lambda^i_\vert(M)$,
hence it is orthogonal to $\Lambda^*_\hor(M)$. This proves \eqref{_d^*_bas_Equation_}. \endproof

\hfill

Now we can prove  \ref{_basic_harmo_Proposition_}.
By \eqref{_d^*_bas_Equation_}, for any basic form $\alpha$ we have 
\begin{equation}\label{_delta_bas_dual_Equation_}
g(\Delta_\bas\alpha, \alpha)= (d\alpha, d\alpha)+ (d^*_\bas \alpha, d^*_\bas \alpha),
\end{equation}
hence a basic form belongs to $\ker \Delta_\bas$ if and only if it is 
closed and orthogonal to all exact basic forms. This gives
an embedding 
\begin{equation}\label{_harmo_to_bas_Equation_}
\ker \Delta_\bas\restrict {\Lambda^*_\bas(M)}\hookrightarrow H^*_\bas(M).
\end{equation}
It remains only to show that the map \eqref{_harmo_to_bas_Equation_} is surjective.
This is where we use the assumption (*) of \ref{_basic_harmo_Proposition_}.

Consider the Hattori decomposition $d=d_0+d_1+\cdots+d_{r+1}$ (see 
Subsection \ref{_Leray_Serre_Subsection_}).
Let $v_1,\ldots, v_r\in F$ be the Killing, transversally Killing vector fields,
postulated in (*), and $\Delta_s$ the ``split Laplacian'',
$\Delta_s:= \{d_1, d_1^*\}- \sum_i \Lie_{v_i}^2$. 
Clearly, $\Delta_s$ is
an elliptic, second order differential operator. 
By definition, $g(e_{v_i} i_{v_i}\alpha, \alpha)= g(i_{v_i}\alpha,i_{v_i}\alpha)$.
Since $v_i$ are Killing, one has $\Lie_{v_i}=-\Lie_{v_i}^*$. Therefore
$\Delta_s$ is self-adjoint and positive, with
\[
g(\Delta_s\alpha, \alpha) = g(d_1\alpha, d_1\alpha) +  
g(d^*_1\alpha, d^*_1\alpha) + \sum_i g(\Lie_{v_i}\alpha,\Lie_{v_i}\alpha). 
\]
We obtain that each $\alpha \in \ker \Delta_s$ satisfies
$i_{v_i} \alpha = \Lie_{v_i}\alpha=0$.

Consider the orthogonal projection map $\Pi_\hor:\; \Lambda^*(M)\arrow \Lambda^*_\hor(M)$.
Since the vector fields $v_i$ are Killing and preserve $F$, 
the projection $\Pi_\hor$ commutes with $\Lie_{v_i}$. It is not hard to see
that $\Pi_\hor$ commutes with $d_1$ and $d_1^*$.
Indeed, $d_1$ is the part of $d$ which maps
$\Lambda^p_\hor(M)\otimes \Lambda^q_\vert(M)$
to $\Lambda^{p+1}_\hor(M)\otimes \Lambda^q_\vert(M)$,
and $d_1^*$ its adjoint. We obtain that 
$\Pi_\hor$ commutes with $\Delta_s$.

Since $\Delta_s$ is a positive, self-adjoint, Fredholm operator,
its eigenvectors are dense in $\Lambda^*(M)$.
Since $[ \Pi_\hor, \Delta_s]=0$, the eigenvectors of 
$\Delta_s$ are dense in $\Lambda^*_\hor(M)$.

Let $G$ be the closure of the Lie group generated by
the action of $e^{\R v_i}$. Since each $e^{\R v_i}$ acts by isometries,
the group $G$ is compact. By construction, $G$ acts on $M$
preserving the foliation $F$, the metric and the transversal
K\"ahler structure, hence it commutes with the Laplacian.
Averaging on $G$, we obtain that
the eigenvectors of $\Delta_s$ are dense in the space
$\Lambda^*_\hor(M)^G= \Lambda^*_\bas(M)$ of all basic forms.

On basic forms, $d=d_1$, hence on $\Lambda^*_\bas(M)$ one has 
$[d, \Delta_s]=0$. Let $\Lambda^*_\bas(M)_\lambda$ be the eigenspace
of $\Delta_s\restrict {\Lambda^*_\bas}(M)$ corresponding to the eigenvalue $\lambda$.
For any closed $\alpha\in \Lambda^*_\bas(M)$
we have $\lambda\alpha = (dd^*+ d^*d)(\alpha)= dd^*\alpha$.
Therefore, any closed form in $\Lambda^*_\bas(M)_\lambda$ is
exact when $\lambda\neq 0$. We have shown that the only eigenspace of $\Delta_\bas$
which contributes to the basic cohomology is 
$\Lambda^*_\bas(M)_0= \ker \Delta_\bas\restrict {\Lambda^*_\bas(M)}$.
This proves \eqref{_harmo_to_bas_Equation_}; we finished the
proof of \ref{_basic_harmo_Proposition_} and \ref{_Transversal_Lefschetz_Theorem_}.
\endproof


\section{Basic cohomology and Hodge theory on Sasakian
  manifolds}
\label{_Hodge_Sasakian_Section_}


\subsection{Cone of a morphism of complexes and cohomology of Sasakian manifolds}
We recall first several notions in homological algebra, see \cite{_Gelfand_Manin_}:

\hfill

\definition
{\bf A complex} $(C_*, d)$ is a collection of vector spaces and homomorphisms
\[ 
\dots\stackrel d \arrow C_i  \stackrel d \arrow C_{i+1} \stackrel d \arrow \cdots
\]
(more generally, a collection of objects in an abelian category)
such that $d^2=0$. {\bf A morphism} of complexes is a collection of maps
$C_i \arrow C_i'$ from the vector spaces of 
a complex $(C_*, d)$ to the vector spaces of $(C'_*, d)$, commuting with the differential.
{\bf The cohomology groups} of a complex $(C_*, d)$  are the groups 
\[ H^i(C_*, d):= \frac{\ker d\restrict C_i}{\im d\restrict{C_{i-1}}}.\]
Clearly, any morphism induces a homomorphism in cohomology.
{\bf An exact sequence of complexes} is a sequence 
\[ 0\arrow A_* \arrow B_* \arrow C_*\arrow 0\] of morphisms 
of complexes such that the corresponding sequences
\[ 0\arrow A_i \arrow B_i \arrow C_i\arrow 0\]
are exact for all $i$.

\hfill

The following claim is very basic.

\hfill

\claim
Let $0\arrow A_* \arrow B_* \arrow C_*\arrow 0$
be an exact sequence of complexes. Then
there is a natural long exact sequence of cohomology
\[
\cdots \rightarrow H^{i-1}(C_*, d_C) \rightarrow H^{i}(A_*, d_A)\rightarrow H^{i}(B_*, d_B)
\rightarrow H^{i}(C_*, d_C) \rightarrow \cdots \quad\text{\endproof}
\]

\definition
Let  $(C_*, d_C)\stackrel \phi \arrow (C'_*, d_{C'})$
be a morphism of complexes. Consider the complex
$C(\phi)_*$, with $C(\phi)_i=C_{i+1}\oplus C'_{i}$
and differential 
\[ d:=d_C+\phi-d_{C'}:\; C_i\oplus C'_{i-1}\arrow 
C_{i+1}\oplus C'_{i}.
\]
or, explicitly:
\[
d(c_i,c'_{i-1})=\left(d_C(c_i), \phi(c_i)-d_{C'}(c_{i-1}')\right).
\]
The complex $(C(\phi)_*, d)$ is called
{\bf the cone of $\phi$}.

\hfill

Denote by $C_*[1]$ the complex $C_*$ shifted by one,
with $C_i[1]=C_{i+1}$. 
The exact sequence of complexes
\[ 
0 \arrow C_*'\arrow C(\phi)_*\arrow C_*[1]\arrow 0
\]
gives the long exact sequence
\begin{multline}\label{_long_exact_from_cone_Equation_}
\cdots \arrow H^{i}(C)\stackrel \phi \arrow H^i(C') 
\arrow  H^i(C(\phi))\arrow\\ \arrow H^{i+1}(C) \stackrel \phi\arrow   H^{i+1}(C') \arrow \cdots
\end{multline}

\proposition\label{_Sasakian_vecr_inva_cone_Proposition_}
Let $Q$ be a Sasakian manifold, and $\vecr$ the Reeb field.
Denote by $\Lambda^*_\vecr(Q)$ the differential
graded algebra of $\Lie_\vecr$-invariant forms,
and let $\Lambda^*_\bas(Q)\subset\Lambda^*_\vecr(Q)$ be the algebra of 
basic forms. Denote by $\omega_0\in \Lambda^*_\bas(Q)$ the transversal 
K\"ahler form (\ref{_omega_0_trans_Kahler_Claim_}), and let 
$L_{\omega_0}:\; \Lambda^*_\bas(Q)\arrow \Lambda^{*+2}_\bas(Q)$
be the multiplication map. Then the complex $\Lambda^*_\vecr(Q)$
is naturally identified with $C(L_{\omega_0})[-1]$,
where $C(L_{\omega_0})$ is the cone of the 
morphism $L_{\omega_0}:\; \Lambda^*_\bas(Q)\arrow \Lambda^{*+2}_\bas(Q)$.

\hfill

\proof
Consider the Hattori decomposition of the differential in $\Lambda^*(Q)$, 
$d= d_0+ d_1 + d_2$ (Subsection \ref{_LS_decompo_Sasakian_Subsection_}).
By \ref{_d_0_expli_Claim_}, $d_0$ vanishes on $\Lambda^*_\vecr(Q)$;
\ref{_d_2_expli_Proposition_} implies that $d_2=L_{\omega_0}i_\vecr$.
Clearly, $\Lambda^*_\vecr(Q)= \Lambda_\bas^*(Q)\oplus \vecr^{\,\flat} \wedge \Lambda_\bas^*(Q)$.
The operator $d_2\Lambda_{\omega_0}i_\vecr$ acts trivially on 
$\Lambda_\bas^*(Q)$ and maps $\vecr^{\,\flat} \wedge \alpha$ to
$L_{\omega_0}(\alpha)$ for any $\alpha \in \Lambda_\bas^*(Q)$.
Therefore, under the natural identification
\[
\Lambda^*_\vecr(Q)= 
\Lambda_\bas^*(Q)\oplus \vecr^{\,\flat} \wedge \Lambda_\bas^*(Q) = \Lambda_\bas^*(Q)\oplus \Lambda_\bas^*(Q)[-1].
\]
By \ref{_d_2_expli_Proposition_}, the differential $d_1+ d_2$ in 
$\Lambda^*_\vecr(Q)$ gives the same differential as in 
$C(L_{\omega_0})[-1]=\Lambda_\bas^*(Q)\oplus \Lambda_\bas^*(Q)[-1]$:
\[
d(\alpha\oplus \vecr^{\,\flat} \wedge\beta)= d\alpha+ L_{\omega_0}(\beta)\oplus (- \vecr^{\,\flat} \wedge d \beta)
\]
for any $\alpha, \beta \in \Lambda_\bas^*(Q)$.
\endproof

\hfill


\subsection{Harmonic forms decomposition on Sasakian manifolds}

We give a new proof of the main result on harmonic forms
on compact Sasakian manifolds (see \cite[Proposition 7.4.13]{_Boyer_Galicki_}, essentially based on a result by Tachibana, \cite{_Tachibana_}).

\hfill

\theorem\label{_Sasakian_cohomo_cone_Theorem_}
Let $Q$ be a $2n+1$-dimensional compact Sasakian manifold,  $R$ the Reeb foliation, 
and $H^*_\bas(Q)$ the corresponding basic cohomology.
Consider the Lefschetz ${\goth {sl}}(2)$-triple
$L_{\omega_0}, \Lambda_{\omega_0}, H_{\omega_0}$ acting
on the basic cohomology (\ref{_Transversal_Lefschetz_Theorem_}).
Then:
\begin{equation*}
	 H^i(Q)= \begin{cases}
	 \ker L_{\omega_0} \restrict {H^i_\bas(Q)}, \qquad \text{for}\ i\geq n,\\[.2in]
	 \dfrac{H^i_\bas(Q)}{\im L_{\omega_0}}, \qquad\qquad \  \text{for}\ i< n.
	 \end{cases}
\end{equation*}

\pstep
The cohomology of the algebra $\Lambda^*_\vecr(Q)$ of $\Lie_\vecr\ $-invariant
forms is equal to the cohomology of $\Lambda^*(Q)$. Indeed,
let $G$ be the closure of the action of $e^{t\vecr}$.
Since $\vecr$ is Killing, $e^{t\vecr}$ acts on $Q$ by isometries,
hence its closure $G$ is compact. Since $G$ is connected
its action on cohomology is trivial. The averaging map
$\Av_G:\; \Lambda^*(Q)\arrow \Lambda^*_\vecr(Q)$
induces an isomorphism on cohomology.

\hfill

{\bf Step 2:} Applying \eqref{_long_exact_from_cone_Equation_},
\ref{_Sasakian_vecr_inva_cone_Proposition_}, and taking into account the isomorphism 
$H^*(\Lambda^*_\vecr(Q))= H^*(Q)$, we obtain the long exact sequence
\begin{multline}\label{_cone_Sasa_equation_}
\cdots \arrow H^{i-2}_\bas(Q)\xrightarrow{L_{\omega_0}}  H^{i}_\bas(Q) 
\arrow  H^{i}(Q)\arrow\\
\arrow H^{i-1}_\bas(Q) \xrightarrow{L_{\omega_0}}   H^{i+1}_\bas(Q) \arrow \cdots
\end{multline}
Since the Lefschetz triple $L_{\omega_0}, \Lambda_{\omega_0}, H_{\omega_0}$
induces an $\goth{sl}(2)$-action on $H^{i}_\bas(Q)$ (\ref{_Transversal_Lefschetz_Theorem_}),
the map $H^{i-2}_\bas(Q)\stackrel {L_{\omega_0}} \arrow H^{i}_\bas(Q)$
is injective for $i\leq n$ and surjective for $i>n$.
Therefore, the long exact sequence \eqref{_cone_Sasa_equation_}
gives the short exact sequences
\[
0\arrow H^{i-2}_\bas(Q)\stackrel {L_{\omega_0}} \arrow H^{i}_\bas(Q) 
\arrow  H^{i}(Q)\arrow 0
\]
for $i\leq n$ and 
\[
0 \arrow H^{i}(Q)\arrow H^{i-1}_\bas(Q) \stackrel {L_{\omega_0}} 
\arrow   H^{i+1}_\bas(Q)\arrow 0
\]
for $i> n$.
\endproof

\hfill

\theorem\label{_harmo_Sasa_decompo_Theorem_}
Let $Q$ be a $2n+1$-dimensional compact Sasakian manifold, and
${\cal H}^i_\bas(Q)$ the space of basic harmonic forms.
Let $\vecr^{\,\flat}$ be the contact form, dual to the Reeb field.
Denote by ${\cal H}^i$ the space of all $i$-forms $\alpha$ on $Q$ 
which satisfy
\begin{description}
\item[for $i \leq n$:]
$\alpha$ is basic harmonic (that is, belongs 
to the kernel of the basic Laplacian) and satisfies
$\Lambda_{\omega_0}(\alpha)=0$; 
\item[for $i > n$:]
 $\alpha = \beta \wedge \vecr^{\,\flat}$
where $\beta$ is basic harmonic and satisfies
$L_{\omega_0}(\beta)=0$. 
\end{description}
Then all elements of ${\cal H}^*$ are harmonic, and, moreover,
all harmonic forms on $Q$ belong to ${\cal H}^*$.

\hfill

\pstep
We prove that all $\gamma\in {\cal H}^*$ are harmonic.
Let $*_\bas$ be the Hodge star operator on basic forms.
Then for any $\gamma \in \Lambda^*_\hor(Q)$ one has
$*(\gamma) = *_\bas(\gamma) \wedge \vecr^{\,\flat}$.
This implies that the  two classes of forms in
\ref{_harmo_Sasa_decompo_Theorem_} are exchanged by
$*$, and it suffices to prove that
all $\alpha\in \Lambda^i(Q)$ which are basic harmonic
and satisfy  $\Lambda_{\omega_0}(\alpha)=0$ 
for $i\leq n$ are harmonic. Such a form $\alpha$ is  closed by \ref{_Transversal_Lefschetz_Theorem_}.
Since $\alpha$ is basic hence $\Lie_\vecr$-invariant, 
it satisfies $d_0^*(\alpha)=0$ (\ref{_d_0_expli_Claim_}).
By \eqref{_delta_bas_dual_Equation_}, a
basic form $\alpha$ is basic harmonic if and only if
$d_1(\alpha)=d_1^*(\alpha)=0$. Finally, $d_2^*=\Lambda_{\omega_0}e_{\vecr}$
(\ref{_d_2_expli_Proposition_}), hence $d_2^*(\alpha)=0$. This implies
$d^*\alpha= (d_0^*+d_1^*+d_2^*)(\alpha)=0$.

\hfill

{\bf Step 2:} 
Now we prove that all harmonic forms on $Q$ belong to ${\cal H}^*$.
By \ref{_Sasakian_cohomo_cone_Theorem_}, 
the dimension of $H^i(Q)$ is equal to the dimension of 
${\cal H}^i$, hence the embedding
${\cal H}^i \arrow H^i(Q)$ constructed in Step 1
is also surjective for all indices $i$.
\endproof


\section{Vaisman manifolds}
\label{_Vaisman_Section_}


In this section we present, without proofs, the necessary background for LCK and Vaisman manifolds. We refer to \cite{do} and to recent papers of ours for details, e.g. \cite{ov_jgp_16}.

\subsection{LCK manifolds}

\definition Let $(M,I,g)$ be a Hermitian manifold of complex dimension $n\geq 2$, with fundamental form $\omega$. It is called {\bf locally conformally K\"ahler (LCK)} if there exists a closed 1-form $\theta$ ({\bf the Lee form}) such that $d\omega=\theta\wedge\omega$. If $\theta$ is exact, $M$ is called {\bf globally conformally K\"ahler (GCK)}. The vector field $\theta^\sharp$ metrically equivalent with $\theta$ is called {\bf the Lee field}.

\hfill

\remark The definition is conformally invariant: if $(M,I,g,\theta)$ is LCK, then $(M,I,e^fg,\theta+df)$ is LCK too. 

\hfill

\remark {(\bf Equivalent definition)}\label{_Equiv_def_LCK_} $(M,I,g)$ is LCK if and only if there exists a cover $\tilde M$ admitting a K\"ahler metric $\tilde g$ such that the deck group  $\Gamma$ acts with holomorphic homotheties w.r.t. $\tilde g$.

\hfill

\remark (i) On a K\"ahler cover as above, the pull-back of the Lee form is exact and the pull-back of the LCK metric is GCK.

(ii) Let $(\tilde M, I, \tilde g)$ be the universal cover of $M$. The equivalent definition in \ref{_Equiv_def_LCK_}  shows the existence of a homothety character $\chi:\pi_1(M)\longrightarrow \R^{>0}$ defined by $\chi(\gamma)=\frac{\gamma^*g}{g}$, for all $\gamma\in\pi_1(M)$. 

\hfill

\definition The rank of $\im(\chi)\subset \R^{>0}$ is called the {\bf LCK rank} of the LCK manifold $M$.

\subsection{Vaisman manifolds}

\definition An LCK manifold $(M,I,g,\theta)$ is called {\bf Vaisman} if the Lee form is parallel w.r.t. the Levi-Civita connection of $g$.

\hfill

\remark\label{_canonical_foliation_properties_} (i) The Lee form of a Vaisman metric is, in particular, co-closed, and hence it is a Gauduchon metric. Therefore, on a compact LCK manifold a Vaisman metric, if it exists, is unique in its conformal class (up to constant multiplier).

(ii) Since $\theta$ is parallel, it has constant norm. The metric $g$ can be rescaled such that $|\theta|=1$. The fundamental form $\omega$ of the Vaisman metric with unit length Lee form satisfies the equality:
\begin{equation}\label{_omega_equation_on_Vaisman_}
d(I\theta)=\omega-\theta\wedge(I\theta).
\end{equation}
In particular, $d(I\theta)$ is positive definite on $\Sigma^\perp$, thus defining a volume form on $\Sigma^\perp$.
(iii) Let $(M,I,g,\theta)$  be a Vaisman manifold,
$\theta^\sharp$ its Lee field. Then $\theta^\sharp$ and
$I\theta^\sharp$ are holomorphic
($\Lie_{\theta^\sharp}I=0$, $\Lie_{I\theta^\sharp}I=0$)
and Killing ($\Lie_{\theta^\sharp}g=0$,
$\Lie_{I\theta^\sharp}g=0$). Therefore, they generate a
foliation $\Sigma$ of real dimension 2 which is complex,
Riemannian and totally geodesic. In particular, $\Sigma$
is transversally K\"ahler (see
\ref{_transvers_Kahler_Definition_}).

(iv) Moreover, $\Sigma$ is canonical in the following
sense: on a compact LCK manifold, the Lee fields of  all
Vaisman metrics are proportional, and hence the foliation
$\Sigma$ is the same for all Vaisman metrics. Therefore,
$\Sigma$ is called {\bf the canonical foliation}. 

(v) Compact complex submanifolds of Vaisman manifolds are Vaisman
(\cite{_Verbitsky:LCHK_}).

\hfill

\ref{_canonical_foliation_properties_}  (iii) admits the following converse which is a powerful criterion for the existence of a Vaisman metric in a conformal class (it was proven in \cite{kor} for compact LCK manifolds, but a careful analysis of the proof shows that it is valid for non-compact manifolds too).

\hfill

\theorem Let $(M,I,g, \theta)$ be an LCK manifold equipped with a 
holomorphic and conformal $\C$-action $\rho$ without fixed points,
which lifts to non-isometric homotheties on 
a K\"ahler cover $\tilde M$. Then $g$
	is conformally equivalent with a Vaisman metric.
	
	\hfill
	
A Vaisman manifold $M$ can have any LCK rank between 1 and $b_1(M)$. However, they always admit deformations to Vaisman structures with LCK rank 1:

\hfill

\theorem (\cite{ov_jgp_16}) \label{_Vaisman_defo_transve_Proposition_}
Let $(M,I, g, \theta)$ be a compact Vaisman manifold,
$\alpha$ a harmonic 1-form such that the deformed 1-form $\theta':= \theta+\alpha$ has rational cohomology class.
Consider the (1,1)-form $\omega':=d_{\theta'}=\theta' \wedge I\theta' - d^c\theta'$ obtained as a
deformation of $\omega=  \theta \wedge I\theta - d^c\theta$ (cf. \eqref{_omega_equation_on_Vaisman_}). Assume
that $\alpha$ is chosen sufficiently small in such a way
that $\omega'$ is positive definite. 
Then the LCK metric associated $\omega'$ is conformally equivalent to a Vaisman metric.

\hfill

\remark It is known that a mapping torus of a compact contact manifold is locally conformally symplectic. Correspondingly, a mapping torus of a compact Sasakian manifold is a Vaisman manifold (\cite{_Vaisman_gd_}). \ref{_Vaisman_defo_transve_Proposition_}  can be used to prove the following converse (the Structure Theorem for compact Vaisman manifolds):

\hfill

\theorem \label{str_vai}(\cite{ov_str, ov_jgp_16}) Every compact Vaisman manifold is biholomorphic 
	to $C(S)/\Z$, where $S$ is Sasakian, 
$\Z= \bigg\langle (x, t) \mapsto (\phi(x), q t)\bigg\rangle$, $q>1$,
 $\phi$ is a Sasakian automorphism of $S$, and $C(S)$ is the Sasakian cone considered as a complex manifold.

\hfill

\remark \label{_prod_decomp_}Since $\theta$ is parallel, the de Rham splitting theorem implies that  a Vaisman manifold is locally the product of $\R$ with a Riemannian manifold which can be shown to be Sasakian. The above Structure Theorem shows that this local decomposition is canonical.

\hfill

\example (i) All linear Hopf manifolds $(\C^n\setminus 0)/\langle A\rangle$, with $A\in \GL(n,\C)$ semi-simple are Vaisman. From  \ref{_canonical_foliation_properties_} (v) it follows that all compact submanifolds of semi-simple linear Hopf manifolds are Vaisman.

(ii) Vaisman compact surfaces are classified by Belgun (see \cite{bel}, \cite{_ovv:surf_}): diagonal Hopf surfaces and elliptic surfaces.

(iii) Non semi-simple Hopf manifolds, Kato manifolds (\cite{_Istrati_Otiman_Pontecorvo_}) and some Oeljeklaus-Toma manifolds are   LCK manifolds which cannot have Vaisman metrics.  


\section{Hodge theory on Vaisman manifolds}
\label{_Hodge_Vaisman_Section_}


\subsection{Basic cohomology of Vaisman manifolds}


Let $M$ be a Vaisman manifold, $\theta^\sharp$ its Lee field,
and $\vecr:= I(\theta^\sharp)$. Using the local decomposition
of $M$ as a product of $\R$ and a Sasakian manifold (\ref{_prod_decomp_}),
$\vecr$ can be identified with the Reeb field on its Sasakian component.
The manifold $M$ is equipped with 2 remarkable foliations,
$\Sigma=\langle \vecr, \theta^\sharp\rangle$ and 
${\cal L}=\langle \theta^\sharp\rangle$. Both of these
foliations satisfy the assumptions of \ref{_basic_harmo_Proposition_} (see also \ref{_canonical_foliation_properties_} (iii)).

Consider the corresponding algebras of basic forms:
\begin{itemize}
\item[(i)] $\Lambda^*_\sas(M)$ -- forms which are basic with respect to ${\cal L}$,
\item[(ii)] $\Lambda^*_\kah(M)$ -- forms which are basic with respect to $\Sigma$.
\end{itemize}
Denote by $H^*_\sas(M)$, $H^*_\kah(M)$ the corresponding basic cohomology 
algebras. 

\hfill

By \ref{_Transversal_Lefschetz_Theorem_}, the algebra $H^*_\kah(M)$
is equipped with the Lefschetz ${\goth {sl}}(2)$-action and the Hodge decomposition
in the same way as for K\"ahler manifolds.

\hfill

The cohomology of Vaisman manifolds can be expressed non-am\-bi\-gu\-ously in terms
of $H^*_\kah(M)$ and the Lefschetz ${\goth {sl}}(2)$-action.
The following theorem is the Vaisman analogue
of \ref{_Sasakian_cohomo_cone_Theorem_}.

\hfill

\theorem\label{_Vaisman_coho_Theorem_}
Let $M$ be a compact Vaisman manifold, $\dim_\R M=2n$, and $\theta$ its
Lee form. Then 
\begin{align}
\label{_coho_Vaisman_via_sasa_Equation_}H^i(M)& \cong H^i_\sas(M) \oplus \theta \wedge H^{i-1}_\sas(M), \\[.1in]
H^i_\sas(M)&=
\begin{cases}
\displaystyle\ker L_{\omega_0} \restrict {H^i_\kah(M)} & \text{for} \ \  i\geq n-1,\\[.2in]
\displaystyle\frac{H^i_\kah(M)}{\im L_{\omega_0}} & \text{for}\ \   i< n-1
\end{cases} \nonumber
\end{align}

\noindent where $\{L_{\omega_0}, \Lambda_{\omega_0}, H_{\omega_0}\}$
is the Lefschetz ${\goth {sl}}(2)$-triple acting on $H^i_\kah(M)$.

\hfill

\pstep We start by proving  \eqref{_coho_Vaisman_via_sasa_Equation_}.
Let $G_{\theta^\sharp}$ be the closure of the
1-parametric group $e^{t\theta^\sharp}$ in the group
$\Iso(M)$ of isometries of $M$. Since
$\Iso(M)$ is compact, the group $G_{\theta^\sharp}$
is also compact. Clearly, averaging on 
$G_{\theta^\sharp}$ does not change the cohomology class
of a form. Therefore, the algebra of
$G_{\theta^\sharp}$-invariant forms has the same
cohomology as $\Lambda^*(M)$. 
Consider the decomposition
\eqref{_horis_vert_decompo_Equation_}
associated with the rank 1 foliation ${\cal L}$:
\begin{equation}\label{_horis_vert_2_decompo_Equation_}
 \Lambda^m(M)= \bigoplus_p \Lambda^p_\hor(M)\otimes
 \Lambda^{m-p}_\vert(M).
\end{equation}
The bundle $\Lambda^{1}_\vert(M)= {\cal L}^*$
is 1-dimensional and generated by $\theta$. Therefore,
\eqref{_horis_vert_2_decompo_Equation_}
gives
\[
\Lambda^m(M)= \Lambda^m_\hor(M)\oplus \theta \wedge \Lambda^{m-1}_\hor(M).
\]
Denote by $\Lambda^m(M)^{G_{\theta^\sharp}}$ the $G_{\theta^\sharp}$-invariant part of $\Lambda^m(M)$. 
Since $\theta$ is $G_{\theta^\sharp}$-invariant, and the
$G_{\theta^\sharp}$-invariant part of $\Lambda^m_\hor(M)$
is identified with $\Lambda^m_\sas(M)$, this gives
\[
\Lambda^m(M)^{G_{\theta^\sharp}}= 
\Lambda^m_\sas(M)\oplus \theta \wedge \Lambda^{m-1}_\sas(M).
\]
Taking cohomology, we obtain \eqref{_coho_Vaisman_via_sasa_Equation_}.

\hfill

{\bf Step 2:} Now we shall prove
$H^i_\sas(M)= \ker L_{\omega_0} \restrict {H^i_\kah(M)}$
for $i\geq n-1$ and $H^i_\sas(M)= \dfrac{H^i_\kah(M)}{\im L_{\omega_0}}$
for $i< n-1$.

Let $\vecr\subset TM$ be the Reeb field defined as above, and
$G_\vecr$ the closure of the 1-parametric group $e^{t \vecr}$.
Then (as in Step 1 and in the proof of \ref{_Sasakian_cohomo_cone_Theorem_}),
the $G_\vecr$-invariant part of $\Lambda^*_\sas(M)$ is written as
\begin{equation}\label{_Sasakian_inva_decompo_Equation_}
\Lambda^i_\sas(M)^{G_\vecr}=\Lambda^i_\kah(M) \oplus \vecr^{\,\flat} \wedge \Lambda^{i-1}_\kah(M)
\end{equation}
with the differential acting (see \ref{_Sasakian_vecr_inva_cone_Proposition_})
as
\begin{equation}\label{_Sasakian_inva_decompo_differe_Equation_}
d(\alpha\oplus \vecr^{\,\flat} \wedge\beta)= d\alpha+ L_{\omega_0}(\beta)\oplus (- \vecr^{\,\flat} \wedge d \beta)
\end{equation}
for any $\alpha, \beta\in \Lambda^*_\kah(M)$.

Using Cartan's formula, we notice that $G_\vecr$ acts trivially
on the cohomology of the complex $(\Lambda^*_\sas(M), d)$, hence the cohomology
of $(\Lambda^*_\sas(M)^{G_\vecr}, d)$ is identified with $H^*_\sas(M)$.
From \eqref{_Sasakian_inva_decompo_Equation_} and
\eqref{_Sasakian_inva_decompo_differe_Equation_}, we obtain that
$\Lambda^i_\sas(M)^{G_\vecr}$ is the cone of the morphism of complexes
$L_{\omega_0}:\; \Lambda^*_\kah(M)[-1] \arrow \Lambda^*_\kah(M)[1]$.
From the long exact sequence
\eqref{_Sasakian_inva_decompo_differe_Equation_}
we obtain a long exact sequence identical to \eqref{_cone_Sasa_equation_}:
\begin{multline}
\cdots \arrow H^{i-2}_\kah(M)\xrightarrow{L_{\omega_0}}  H^{i}_\kah(M) 
\arrow  H^{i}_\sas(M)\arrow \\
\arrow H^{i-1}_\kah(M) \xrightarrow{L_{\omega_0}}   H^{i+1}_\kah(M) 
\arrow\cdots
\end{multline}
Since $L_{\omega_0}:\; H^{i-2}_\kah(M)\arrow  H^{i}_\kah(M)$ is injective for $i\leq n-1$
and surjective for $i >n-1$, this long exact sequence breaks into short
exact sequences of the form
\[
0\arrow H^{i-2}_\kah(M)\stackrel {L_{\omega_0}} \arrow H^{i}_\kah(M) 
\arrow  H^{i}_\sas(M)\arrow 0
\]
for $i\leq n-1$ and 
\[
0 \arrow H^{i}(M)_\sas \arrow H^{i-1}_\kah(M) \stackrel {L_{\omega_0}} \arrow   H^{i+1}_\kah(M)\
\arrow 0
\]
for $i> n$. We finished the proof of 
\ref{_Vaisman_coho_Theorem_}. \endproof

\subsection{Harmonic forms on Vaisman manifolds}

It turns out that (just like it happens in the Sasakian case)
the cohomology decomposition obtained in \ref{_Vaisman_coho_Theorem_}
gives a harmonic form decomposition. Together with the Hodge decomposition of the 
basic cohomology of transversally K\"ahler structure this allows us to represent
certain cohomology classes by forms of a given Hodge type. 
This theorem was obtained in \cite{_Kashiwada_Sato_}, see also \cite{_Vaisman_gd_} and \cite{tsu}.

\hfill

Recall that on a compact Vaisman manifold $(M,I,g,\theta)$ with fundamental form $\omega$, the canonical foliation $\Sigma$ is transversally K\"ahler (\ref{_canonical_foliation_properties_} (iii)). We denoted $\vecr$ the vector field $g$-equivalent with $I\theta$, $\Lambda^*_\kah(M)$ the space of 
basic forms with respect to $\Sigma$, and 
$\Delta_\kah:\; \Lambda^*_\kah(M)\arrow \Lambda^*_\kah(M)$
the transversal Laplacian (\ref{_basic_harmo_Proposition_}). From \ref{_Transversal_Lefschetz_Theorem_}
we obtain that the space ${\cal H}^i_\kah(M)$ 
of basic harmonic forms is equipped with the
Lefschetz ${\goth {sl}}(2)$-action by the operators 
$L_{\omega_0}, \Lambda_{\omega_0}, H_{\omega_0}$. The main result of this section is:

\hfill

\theorem\label{_Vaisman_harmonic_forms_Theorem_}
Let $(M,I,g,\theta)$ be a compact Vaisman manifold of complex dimension $n$, with fundamental form $\omega$, and canonical foliation $\Sigma$. 
Denote by ${\cal H}^i$ the space of all  basic $i$-forms 
$\alpha\in \Lambda^*_\kah(M)$
which satisfy:
\begin{description}
\item[\hspace{.1in}for $i \leq n$:]
$\alpha$ is basic harmonic (i.e.  $\Delta_\kah(\alpha)=0$) and satisfies
$\Lambda_{\omega_0}(\alpha)=0$; 
\item[\hspace{.1in}for $i > n$:]
 $\alpha = \beta \wedge I\theta$
where $\beta$ is basic harmonic and satisfies
$L_{\omega_0}(\beta)=0$. 
\end{description}
Then all elements of ${\cal H}^*\oplus \theta\wedge {\cal H}^*$ are harmonic and, moreover,
all harmonic forms on $M$ belong to ${\cal H}^*\oplus \theta\wedge {\cal H}^*$.

\hfill

\pstep
This statement is similar to \ref{_harmo_Sasa_decompo_Theorem_},
and the proof is essentially the same. We start by
proving that all $\gamma\in {\cal H}^*\oplus\theta\wedge {\cal H}^*$ are harmonic.
Notice that a product of a parallel form $\rho$ and a harmonic
form is again harmonic (see e.g. \cite[Proposition 2.7]{_Verbitsky:G2_forma_}).
This result is proven in the same way as \ref{_kah_susy_Theorem_},
by considering the Lie superalgebra generated
by $d, d^*$, and multiplication by $\rho$.
Since $\theta$ is parallel, it suffices only to show
that all elements in ${\cal H}^*$ are harmonic:
\begin{equation}\label{_H^*_harmo_Equation_}
{\cal H}^*\subset \ker \Delta.
\end{equation}
Consider the foliation ${\cal L}=\langle \theta^\sharp\rangle$.
Since $M$ is locally a product of a Sasakian manifold
and a line (\ref{_prod_decomp_}), 
all ${\cal L}$-basic, transversally harmonic forms are harmonic.
Therefore, \eqref{_H^*_harmo_Equation_} would follow if we prove
\begin{equation}\label{_H^*_transve_sasa_harmo_Equation_}
{\cal H}^*\subset \ker \Delta_\sas.
\end{equation}
where $\Delta_\sas$ is the transversal Laplacian of 
the foliation ${\cal L}$.

It would suffice to prove \eqref{_H^*_transve_sasa_harmo_Equation_}
locally in $M$. However, locally  ${\cal L}$
admits a leaf space $Q$,
which is Sasakian, and we can regard elements
from ${\cal H}^*$ as forms on $Q$ and $\Delta_\sas$ as the
usual Laplacian on $\Lambda^*(Q)$. 
Then \eqref{_H^*_transve_sasa_harmo_Equation_} follows from
\ref{_harmo_Sasa_decompo_Theorem_}, 
step 1.\footnote{The statement of \ref{_harmo_Sasa_decompo_Theorem_}
is global, however, the proof of \ref{_harmo_Sasa_decompo_Theorem_}, 
step 1 is local, and this statement is essentially identical to 
\eqref{_H^*_transve_sasa_harmo_Equation_}.}

\hfill

{\bf Step 2:} We have shown that all elements of 
${\cal H}^*\oplus \theta\wedge {\cal H}^*$ are harmonic;
this gives a natural linear map
\[
\Psi:\; {\cal H}^*\oplus \theta\wedge {\cal H}^*\arrow H^*(M).
\]
To show that all harmonic form are obtained this way,
it would suffice to prove that $\Psi$ is surjective.
This follows from \ref{_Vaisman_coho_Theorem_}
by dimension count.

By  \ref{_Vaisman_coho_Theorem_}, the cohomology of $M$ is isomorphic to
$H^i_\sas(M) \oplus \theta \wedge H^{i-1}_\sas(M)$.
All forms in ${\cal H}^*$  are closed and belong to $\Lambda^*_\sas(M)$,
which gives a  map $\Psi_\sas:\; {\cal H}^*\arrow H^i_\sas(M)$. 
To prove that $\Psi$ is surjective, it remains to show that
$\Psi_\sas$ is surjective.

However, the dimension of $H^i_\sas(M)$ is equal to
$\dim \ker L_{\omega_0}\restrict{H^i_\kah(M)}$ for $i \geq n-1$ and to 
$\dim \dfrac{H^i_\kah(M)}{\im L_{\omega_0}}$ for $i < n-1$.
The space ${\cal H}^i$ has the same dimension by 
the transversal Hodge decomposition
(\ref{_Transversal_Lefschetz_Theorem_}).
This finishes the proof of 
\ref{_Vaisman_harmonic_forms_Theorem_}.
\endproof


\section{Supersymmetry on Sasakian manifolds}
\label{_SUSY_Sasakian_Section_}


Let $Q$ be a Sasakian manifold.
In this section  we describe the Lie superalgebra $\goth q\subset \End(\Lambda^*(Q))$ 
reminiscent of the supersymmetry algebra of a K\"ahler manifold 
(\ref{_kah_susy_Theorem_}). Unlike the K\"ahler supersymmetry algebra,
the algebra $\goth q$ is infinitely-dimensional; however, it has a simple
and compact description, independent of the choice of $Q$.

When this project was started, we expected to use 
\ref{_susy_Sasa_Theorem_} to give a more conceptual proof of 
the classical results on Hodge decomposition of the
cohomology of Sasakian and Vaisman manifolds (\ref{_harmo_Sasa_decompo_Theorem_}, \ref{_Vaisman_harmonic_forms_Theorem_}). However,
the Lie superalgebra $\goth q$ which we obtained 
in the end does not contain the 
de Rham Laplacian operator (the ``de Rham Laplacian'' is
the usual Laplacian operator defined on differential
forms on a Riemannian manifold). The de Rham operator
on a Sasakian manifold is decomposed onto its
Hattori components as $d_0+ d_1+ d_2$ (Section
\ref{_LS_decompo_Sasakian_Subsection_}). The operators 
$d_0, d_1$ are elements of $\goth q$, but $d_2$ is in fact an
element of its universal enveloping algebra $U_{\goth q}$.

Any attempt to
add $d_2$ or the de Rham Laplacian to $\goth q$ lead
to a large subalgebra of $U_{\goth q}$ which
is complicated and very difficult to control.

The current version of the proof of 
\ref{_Vaisman_harmonic_forms_Theorem_} is independent
from \ref{_susy_Sasa_Theorem_}. 

\hfill

\remark Some of the relations we derive can be found in P.A.-Nagy's doctoral thesis, in a more general setting, for a compact Riemannian manifold endowed with a unitary vector field (see \cite{_Nagy_}).

\hfill


We first recall the notations. Let $Q$ be a Sasakian manifold, $\vecr$ its Reeb field,
and $R\subset TM$ the corresponding rank 1 distribution. 
Let 
$$\Lambda^*(Q):= \bigoplus\Lambda^{p,q}_\hor(Q)\times
\Lambda^m_\vert(Q)$$
 be the corresponding decomposition
(Hodge and Hattori) of the de Rham algebra, and 
$d=d_0+d_1+d_2$ the Hattori differentials (see Subsection
\ref{_LS_decompo_Sasakian_Subsection_}).
Denote by $W$ the Weil operator acting as multiplication
by  $\1(p-q)$ on
$\bigoplus\Lambda^{p,q}_\hor(Q)\times\Lambda^m_\vert(Q)$. 
Let $d_1, d_1^*, d_1^c$ and $(d_1^c)^*$ be the
differentials defined in \ref{_Hodge_decompo_d_1_Sasa_Claim_}.
Let $L_{\omega_0}, \Lambda_{\omega_0}:= L_{\omega_0}^*, H_{\omega_0}:= [L_{\omega_0}, \Lambda_{\omega_0}]$
be the Lefschetz ${\goth {sl}}(2)$-triple associated with the
transversally K\"ahler form $\omega_0$.
For an operator $A\subset \End(Q)$, 
commuting with $\Lie_\vecr$, we denote by
$A(k)$ the composition $A\circ (\Lie_\vecr)^k$.
Let $i_\vecr$ be the
contraction with $\vecr$ and $e_\vecr$ the dual operator (above denoted $e_{\vecr^{\, \flat}}$).

Denote by $\goth q$ the Lie superalgebra generated
by all operators 
\[ L_{\omega_0}(i),\ 
\Lambda_{\omega_0}(i),\ H_{\omega_0}(i), \ d_1(i),\  e_\vecr(i),\ 
i_\vecr(i),\  \Id(i), W(i)\ \text{for all $i\in \Z^{\geq 0}$}.
\] 
Since the vector field $\vecr$ acts by Sasakian
isometries, $\goth q$ commutes with $\Lie_\vecr$;
in other words,  $\Lie_\vecr$ is central in $\goth q$.
We consider $\goth q$ as an $\R[t]$-module, with
$t$ mapping $A(i)$ to $A(i+1)$.
Then the Lie superalgebra $\goth q$ is a free 
$\R[t]$-module of rank $(6|6)$ (with 6 even and 6 odd
  generators) over $\R[t]$. 
Its even generators are 
\[ L_{\omega_0},\ 
\Lambda_{\omega_0},\  H_{\omega_0},\  W,\  \Delta_1:= \{d_1,
d_1^*\},\  \Id,
\]
and the odd generators are 
$$d_1, \ d_1^*,\  d_1^c,\  (d_1^c)^*,\  e_\vecr,\  i_\vecr.$$
Notice that $d_0=e_\vecr(1)$ and $d_0^*=i_\vecr(1)$ (\ref{_d_0_expli_Claim_}).

The main result of this section is:

\hfill

\theorem \label{_susy_Sasa_Theorem_}The only non-zero commutator relations in $\goth q$
can be written as follows.
\begin{description}
\item[(i)] The Lefschetz's ${\goth {sl}}(2)$-action: the even elements
  $L_{\omega_0},\Lambda_{\omega_0},H_{\omega_0}$ 
satisfy the usual ${\goth {sl}}(2)$-relations (\ref{_kah_susy_Theorem_}) 
and commute with 
\[
W, \Delta_1, \Delta_0, d_0, d_0^*.
\]
 The operator
$H_{\omega_0}$ acts as multiplication by 
$p-n$ on $\Lambda^{p}_\hor(Q)\times\Lambda^m_\vert(Q)$,
where $\dim_\R Q= 2n+1$.

\item[(ii)] The Weil operator satisfies 
\[
[W, d]=d_1^c, \ [W, d^c_1]= - d_1,\ 
  [W, d_1^*]=-(d_1^c)^*,\  [W, (d^c_1)^*]= d_1.
\]
Also, $W$ commutes with the rest of the generators of $\goth  q$.

\item[(iii)] The differentials $d_1, d_1^*, d_1^c, (d_1^c)^*$ have non-zero square:
\[ 
 \{d_1, d_1\}=\{d_1^c, d_1^c\}=-L_{\omega_0}(1), \ \ 
\{d_1^*, d_1^*\}=\{(d_1^c)^*, (d_1^c)^*\}=\Lambda_{\omega_0}(1)
\]
Moreover, $\{d_1, d_1^c\}=\{d_1^*, (d_1^c)^*\}=0$.
\item[(iv)] 
The usual Kodaira relations still hold:
\begin{equation}\label{_Kodaira_rel_Sasakian_Equation_}
\begin{split}
  [\Lambda_{\omega_0}, d_1]& = (d_1^c)^*,\ \ \  
\ \ [ L_{\omega_0}, d_1^*] = - d_1^c,\\
[\Lambda_{\omega_0}, d_1^c]& = - d_1^*,
\ \ \ \ \ [ L_{\omega_0}, (d^c_1)^*] = d_1.
\end{split}
\end{equation}

\item[(v)] Unlike it happens in the K\"ahler case,
  the differentials $d_1^*, d_1^c$, etc. do not (super-)commute:
\begin{equation}\label{_commutators_v_}
\{d_1^*, d_1^c\}=\{d_1,
(d_1^c)^*\}= -\frac 1 2 H_{\omega_0}(1).
\end{equation}
\item[(vi)] The only non-zero commutator  between
$e_\vecr, i_\vecr, \Id$ is $\{e_\vecr, i_\vecr\}=\Id$.
These elements commute with the rest of $\goth q$.

\item[(vii)] The Laplacian $\Delta_1:= \{d_1, d_1^*\}$ 
satisfies $\Delta_1= \{ d_1^*, (d_1^c)^*\}$
and commutes with all
even generators in $\goth q$ and with $e_\vecr, i_\vecr$.
Its commutators with the other 4 odd generators are expressed as
follows
\begin{equation}\label{_Delta_1_Equation_}
\begin{split}
 \{d_1,\Delta_1\}&= -\frac 12 d_1^c(1), \ \ \ \ \ \ \ \ \  \{d_1^c,\Delta_1\}= \frac 12
 d_1(1),\\
 \{d_1^*,\Delta_1\}&= - \frac 12 (d_1^c)^*(1),  
 \ \  \{(d_1^c)^*,\Delta_1\}= \frac 12 (d_1)^*(1).
 \end{split}
\end{equation}
\end{description}

{\bf Proof of \ref{_susy_Sasa_Theorem_} (i):}
The ${\goth {sl}}(2)$-relations and the expression for 
$H_{\omega_0}$ are proven in Subsection
\ref{_susy_Kah_Subsection_}; the proof in the Sasakian
case is literally the same. Also, from the definition
it is clear that
$L_{\omega_0},\Lambda_{\omega_0},H_{\omega_0}$ 
commute with $W$, $i_\vecr$ and $e_\vecr$.
Since $d_0=e_\vecr(1)$ and $d_0^*=i_\vecr(1)$ 
(\ref{_d_0_expli_Claim_}), these operators
commute with the $\goth{sl}(2)$-action
and satisfy $\{d_0, d_0^*\}=\Lie^2_\vecr$.
We postpone the commutator relation for $\Delta_1$ until
we proved \ref{_susy_Sasa_Theorem_} (iv).

\hfill

{\bf Proof of \ref{_susy_Sasa_Theorem_} (ii):}
Same as \ref{_kah_susy_Theorem_}, part 2.

\hfill

{\bf Proof of \ref{_susy_Sasa_Theorem_} (iii):}
Start from $\{d_1, d_1\}=-L_{\omega_0}(1)$.

Since $d^2=0$, one has $d_0^2=d_2^2=0$ and
$\{d_0, d_2\}=-\{d_1, d_1\}$. However,
$d_0= e_{\vecr}(1)$ and $d_2=L_{\omega_0} i_\vecr$,
hence $\{d_0, d_2\}=\{e_{\vecr}, i_\vecr\}
L_{\omega_0}(1)=L_{\omega_0} (1)$.
The squares of the rest of the differentials $d_1^c$, etc.,
are obtained by duality and complex conjugation.

To show that $\{d_1, d_1^c\}=0$, we use
$d_1^c=[W, d_1]$. By \ref{_d^2neq0_super_Claim_},
\[
\{d_1, d_1^c\}= \{ d_1, \{d_1, W \}\}= \frac 1 2 \{\{ d_1, d_1\}, W\}.
\]
Since $\{ d_1, d_1\}=-L_{\omega_0}(1)$ and $W$ commutes
with $L_{\omega_0}$, this gives $\{d_1, d_1^c\}=0$. The
equation $\{d_1^*, (d_1^c)^*\}=0$ is obtained by duality.
We finished \ref{_susy_Sasa_Theorem_} (iii).

\hfill

{\bf Proof of \ref{_susy_Sasa_Theorem_} (iv), the K\"ahler-Kodaira relations:}
Again, it suffices to prove 
$ [L_{\omega_0}, d_1^*] = - d_1^c$,
the rest is obtained by duality and complex conjugation.
We prove it applying the same argument as used in
\ref{_kah_susy_Theorem_}.
As an operator on $\Lambda^*(Q)$, the 
commutator  $[L_{\omega_0}, d_1^*]$ has first order,
because $L_{\omega_0}$ is zero order, and $d_1^*$ is order
2 (\ref{_diff_ope_commutator_Claim_}).

The differential operators $[L_{\omega_0}, d_1^*]$
and $- d_1^c$ are equal on functions for the same
reasons as in \ref{_kah_susy_Theorem_}.

Clearly, $d_1^c(C^\infty (Q))$ generates
$\Lambda^1_\hor(Q)$.  Therefore, to prove  
$[L_{\omega_0}, d_1^*]=- d_1^c$ on 
$\Lambda^1_\hor(Q)$, we need only to show that
\begin{equation}\label{_commu_L_d_1^*_square_Equation_}
([L_{\omega_0}, d_1^*])^2= (- d_1^c)^2.
\end{equation}
This is implied by the graded Jacobi identity, applied as follows.
First, we notice that $[\Lambda_{\omega_0}, d_1^*]=0$,
because $d_1(\omega_0)=0$. Therefore, the $\goth{sl}(2)$-representation
generated by the Lefschetz triple 
$\langle
L_{\omega_0},\Lambda_{\omega_0},H_{\omega_0}\rangle$
from $d_1^*$ has weight 1, and
$[L_{\omega_0},[L_{\omega_0}, d_1^*]]=0$.
Applying the graded Jacobi identity, and using
$[L_{\omega_0},[L_{\omega_0}, d_1^*]]=0$,
we obtain

\begin{equation*}
\begin{split}
\{[L_{\omega_0}, d_1^*], [L_{\omega_0}, d_1^*]\}&=
[L_{\omega_0},\{d_1^*, [L_{\omega_0}, d_1^*]\}]\\[.1in] 
&=\frac 1 2[L_{\omega_0}(1), [L_{\omega_0}(1), \{d_1^*,d_1^*\}]]\quad \text{by \ref{_d^2neq0_super_Claim_}}\\
&=-\frac 1 2 [L_{\omega_0}(1), [L_{\omega_0}(1), \Lambda_{\omega_0}]]\quad \text{by \ref{_susy_Sasa_Theorem_} (iii)}\\
&=-L_{\omega_0}(1).
\end{split}
\end{equation*}

However, $(- d_1^c)^2=-L_{\omega_0}(1)$ by
\ref{_susy_Sasa_Theorem_} (iii), which proves
\eqref{_commu_L_d_1^*_square_Equation_}.

This implies that  
$ [L_{\omega_0}, d_1^*] = - d_1^c$
on $\Lambda^1_\hor(Q)$ .
The space $\Lambda^1(Q)$ is generated over
$C^\infty (Q)$ by $\Lambda^1_\hor(Q)$ and 
$V= \langle \vecr^{\,\flat} \rangle$.
Applying \ref{_d^2=0_dete_Corollary_}, we obtain that
$ [L_{\omega_0}, d_1^*] = - d_1^c$ if 
$ [L_{\omega_0}, d_1^*](\vecr^{\,\flat})=-d_1^c(\vecr^{\,\flat})$.
However, $d(\vecr^{\,\flat})=\omega_0$, hence
$d_1^c(\vecr^{\,\flat})=0$. Then, 
$ [L_{\omega_0}, d_1^*](\vecr^{\,\flat})=
d_1^*(\vecr^{\,\flat}\wedge \omega_0)$.
It is not hard to see that 
$*(\vecr^{\,\flat}\wedge \omega_0)=
\frac{1}{(n-1)!}\omega_0^{n-1}$, where $\dim_\R Q=2n+1.$
Therefore, 
\[ d^*_1(\vecr^{\,\flat}\wedge \omega_0) =\pm
*d_1*(\omega_0\wedge \vecr^{\,\flat})=
\pm *d_1 \frac{1}{(n-1)!}\omega_0^{n-1}=0
\]
because $d\omega_0^{n-1}=0$.
We proved that $ [L_{\omega_0},
  d_1^*](\vecr^{\,\flat})=-d_1^c(\vecr^{\,\flat})$
and finished the proof of the K\"ahler-Kodaira relations.

\hfill

{\bf Proof of \ref{_susy_Sasa_Theorem_} (v), the commutators of $d_1^*, d^c$:}
The commutator $\{d_1^*, d^c\}$ is obtained from the
K\"ahler-Kodaira relations. Indeed,
$\{d_1^*, d^c\}=\{d_1^*, \{d_1^*, L_{\omega_0}\}\}$.
Then \ref{_d^2neq0_super_Claim_} gives
\[ \{d_1^*, \{d_1^*, L_{\omega_0}\}= \frac 1 2 \{\{d_1^*,
    d_1^*\}, L_{\omega_0}\}= \frac 1 2 [\Lambda_{\omega_0}(1),
      L_{\omega_0}]=- \frac 1 2 H_{\omega_0}.
\]
The relation $\{d_1,
(d_1^c)^*\}=-\frac 1 2 H_{\omega_0}(1)$ is dual to
$\{d_1^*,d_1^c\}=-\frac 1 2 H_{\omega_0}(1)$.

\hfill

{\bf Proof of \ref{_susy_Sasa_Theorem_} (vi), the commutators of $e_\vecr$ and $i_\vecr$:}
The equation $\{e_\vecr, i_\vecr\}=1$ is standard. Vanishing of the commutators
between $e_\vecr$, $i_\vecr$ and  $L_{\omega_0},
\Lambda_{\omega_0},H_{\omega_0}, W$ is standard linear algebra.
The only commutators for which we have to prove the vanishing
is between $e_\vecr$, $i_\vecr$ and $d_1, d_1^*, d_1^c$, $(d_1^c)^*$.
Using duality and complex conjugation, we reduce the vanishing
of these commutators to only two of them:
$\{e_\vecr, d_1\}=0$ and $\{e_\vecr, d_1^*\}=0$.
As $d^2=0$ and $\{d_0, d_1\}$ is the grading 1 part of
$d^2$, one has $\{d_0, d_1\}=0$. Since $d_0= e_\vecr(1)$,
this also implies $\{e_\vecr, d_1\}=0$. Twisting with $I$, we obtain
$\{e_\vecr, d_1^c\}=0$. Applying the graded Jacobi identity to
$\{e_\vecr, d_1^*\}= -\{ e_\vecr, \{\Lambda_{\omega_0}, d_1^c\}\}$
(\ref{_susy_Sasa_Theorem_} (iv)) and using
 $\{e_\vecr, L_{\omega_0}\}=0$, we obtain
\[
\{e_\vecr, d_1^*\}= -\{ e_\vecr, \{\Lambda_{\omega_0}, d_1^c\}\}=
\{\{e_\vecr, \Lambda_{\omega_0}\},d_1^c\}+ \{\Lambda_{\omega_0}, \{e_\vecr, d_1^c\}\}=0.
\]
This finishes the proof of \ref{_susy_Sasa_Theorem_} (vi).

\hfill

{\bf Proof of   \ref{_susy_Sasa_Theorem_} (vii).}
The equation
\[ \{d_1, d_1^*\}=\{(d_1^c), (d_1^c)^*\}
\]
follows from $\{d_1, d_1^c\}=0$ (\ref{_susy_Sasa_Theorem_} (iii))
because 
\[
0=\{\Lambda, \{d_1, d_1^c\}\}= \{\{\Lambda,d_1\}, d_1^c\}+ \{d_1, \{\Lambda,d_1^c\}\}=
\{(d_1^c)^*,  d_1^c\}-\{d_1, d_1^*\}.
\]
This implies, in particular, that $[W,\{d_1, d_1^*\}]=0$.
The commutators between the Lefschetz operators and $\Delta_1$
follow from the K\"ahler-Kodaira relations:
\[
\{L_{\omega_0}, \{d_1, d_1^*\}\}= \{\{L_{\omega_0}, d_1\}, d_1^*\}
+ \{d_1, \{L_{\omega_0},d_1^*\}\}=-\{d_1,d_1^c\}=0.
\]
We proved that $\Delta_1$ commutes with the even part of $\goth q$.
By duality and complex conjugation, to prove \eqref{_Delta_1_Equation_}
it would suffice to prove only one of these relations, say,
$\{d_1,\Delta_1\}= -\frac 12 d_1^c(1)$. This equation 
follows from \eqref{_Kodaira_rel_Sasakian_Equation_},  \eqref{_commutators_v_} (i.e. (iv) and (v) of this theorem) and \ref{_d^2neq0_super_Claim_}:
\[
\{d_1,\{d_1, d_1^*\}\}=\frac 1 2 \{\{d_1, d_1\}, d_1^*\}=
-\frac 1 2 \{L_{\omega_0}(1), d_1^*\}=-\frac 1 2 (d_1^c)^*.
\]
We finished the proof of \ref{_susy_Sasa_Theorem_}. \endproof

\hfill

{\bf Acknowledgment:} L.O. thanks IMPA (Rio de Janeiro)
and HSE (Mos\-cow) for financial support and excellent research environment
during the preparation of this paper. Both authors thank
P.-A. Nagy for useful discussions. We are grateful to
Richard Eager for the reference to \cite{_Schmude_}
and to Nikita Klemyatin for pointing out some errors in a first version of the paper. Many thanks to the referee and the editor for very useful comments.

\hfill

{\small

\noindent {\sc Liviu Ornea\\
University of Bucharest, Faculty of Mathematics and Informatics, \\14
Academiei str., 70109 Bucharest, Romania}, and:\\
{\sc Institute of Mathematics ``Simion Stoilow" of the Romanian
Academy,\\
21, Calea Grivitei Str.
010702-Bucharest, Romania\\
\tt lornea@fmi.unibuc.ro,   liviu.ornea@imar.ro}

\hfill

\noindent {\sc Misha Verbitsky\\
{\sc Instituto Nacional de Matem\'atica Pura e
              Aplicada (IMPA) \\ Estrada Dona Castorina, 110\\
Jardim Bot\^anico, CEP 22460-320\\
Rio de Janeiro, RJ - Brasil }\\
also:\\
Laboratory of Algebraic Geometry, \\
Faculty of Mathematics, National Research University 
Higher School of Economics,
6 Usacheva Str. Moscow, Russia}\\
\tt verbit@verbit.ru, verbit@impa.br }

\end{document}